\documentclass[a4paper]{article}
\usepackage[top=25truemm,bottom=25truemm,left=20truemm,right=20truemm]{geometry}
\usepackage{amsmath,amsthm,amsfonts,amssymb,stmaryrd,indentfirst,tabularray,makecell,hhline, comment}
\usepackage[mathscr]{euscript}
\usepackage[shortlabels]{enumitem}
\usepackage[labelfont=rm]{caption}
\usepackage{titlesec}
\usepackage{tikz-cd}
\usepackage{url}
\usepackage[backref]{hyperref}
\hypersetup{
    colorlinks = true,
    linkcolor = [rgb]{0.3,0.55,0.25},
    citecolor = [rgb]{0.164,0.39,0.77},
    urlcolor = [rgb]{0.7,0.2,0.4}
}

\DeclareMathOperator{\id}{id}
\DeclareMathOperator{\Ker}{Ker}

\DeclareMathOperator{\Der}{Der}
\DeclareMathOperator{\DDer}{DDer}
\DeclareMathOperator{\Hom}{Hom}

\DeclareMathOperator{\End}{End}

\DeclareMathOperator{\Tr}{Tr}
\DeclareMathOperator{\Div}{Div}
\DeclareMathOperator{\TDiv}{TDiv}
\DeclareMathOperator{\sdiv}{div} % stands for single div

\DeclareMathOperator{\ad}{ad}

\DeclareMathOperator{\trip}{trip}
\DeclareMathOperator{\Ob}{Ob}

\newcommand{\pair}[2]{\langle#1,#2\rangle}
\newcommand{\op}[0]{^\mathrm{op}}
\newcommand{\argdot}[0]{\,\cdot\,}

\theoremstyle{plain}
\newtheorem{theorem}{Theorem}
\newtheorem*{theorem*}{Theorem}
\newtheorem{proposition}[theorem]{Proposition}
\newtheorem{deflem}[theorem]{Definition-Lemma}
\newtheorem{lemma}[theorem]{Lemma}
\newtheorem{corollary}[theorem]{Corollary}

\theoremstyle{definition}
\newtheorem{definition}[theorem]{Definition}
\newtheorem{remark}[theorem]{Remark}
\newtheorem{example}[theorem]{Example}

\numberwithin{theorem}{section}
\everymath{\displaystyle}
\allowdisplaybreaks[1]

\titleformat{\section}{\large\scshape}{\thesection.}{3pt}{}
\titleformat{\subsection}{\scshape}{\thesubsection}{3pt}{}

\usepackage{titling}  
\pretitle{
	\begin{center}\Large\bfseries
}
\posttitle{
	\end{center}\ \\[-20pt]
}
\preauthor{
	\begin{center}\large
}
\postauthor{
	\end{center}\ \\[-60pt]
}

\title{Modular Vector Fields in Non-commutative Geometry}
\author{Toyo TANIGUCHI \thanks{Graduate School of Mathematical Sciences, The University of Tokyo. 3-8-1, Komaba, Meguro-ku, Tokyo, 153-8914, Japan. E-mail: \texttt{toyo(at)ms.u-tokyo.ac.jp}}}
\date{}

\begin{document}
\maketitle

\begin{abstract}
\noindent We construct a non-commutative analogue of the modular vector field on a Poisson manifold for a given pair of a double bracket and a connection on a space of $1$-forms. The key ingredient, the triple divergence map, is directly constructed from a connection on a linear category to deal with multiple base points. As an application, we give an algebraic description of the framed, groupoid version of Turaev's loop operation $\mu$ similar to the one obtained by Alekseev--Kawazumi--Kuno--Naef and the author. \\
\end{abstract}


\section{Introduction}

The Turaev cobracket is a \textit{loop operation} introduced in \cite{turaev} together with its quantisation and gives a Lie bialgebra structure on the vector space spanned by homotopy classes of free loops on a connected compact oriented surface. This Lie bialgebra is a topological counterpart of the Batalin-Vilkovisky structure on the moduli space of flat principal super Lie group bundles; the Lie algebra part is classically known by Goldman \cite{goldman}, and the extension to the whole BV structure was done by Alekseev, Naef, Pulmann, and \v{S}evera \cite{anps}.

An application of the Turaev cobracket in the theory of mapping class groups, and in particular the Johnson homomorphisms and the Enomoto--Satoh trace, was explored with the help of the deep connection to the Kashiwara--Vergne problem in the paper \cite{akkn} by Alekseev, Kawazumi, Kuno and Naef. Along the way, they factored the Turaev cobracket in two components: the ``Hamiltonian flow'' $\sigma$ and the divergence map $\Div^\mathcal{C}$ associated with the free generating system $\mathcal{C}$ of the fundamental group of such a surface with non-empty boundary. In a previous paper \cite{toyo} of the author, the latter was reformulated in terms of non-commutative geometry: associated with a connection $\nabla$ on a finitely generated projective $A$-bimodule $M$ with a derivation action, the divergence map $\Div^\nabla$ is defined to yield a corresponding factorisation in the case of closed surfaces.

In this paper, we will extend the construction to give a similar description of the framed version $\mu^\mathsf{fr}$ of Turaev's loop operation $\mu$ by introducing an analogue of \textit{modular vector fields} in Poisson geometry, which is a characteristic quantity of a Poisson manifold. It enables us to discuss the skew-symmetry of the map $\delta^{\sigma,\nabla}:=\Div^\nabla\circ\,\sigma$, which corresponds to the vanishing of the divergence in the usual Poisson geometry. This skew-symmetry is obvious in some cases since we know that it coincides with the Turaev cobracket for a special connection $\nabla\!_\mathcal{C}$ (see Theorem \ref{thm:cob}), but not at all from the definition of $\delta^{\sigma,\nabla}$ itself.

\begin{comment}
More specifically, denoting the cyclic quotient of an algebra $A$ by $|A|$, we will construct two maps
\begin{align*}
	\phi_{\Pi,\nabla}\colon A \to |\mathrm{trip}(A)|\cong A\otimes |A| \oplus |A|\otimes A \;\textrm{ and }\; \mathbf{m}_{\Pi,\nabla,\Theta} \colon A\to A\otimes |A|
\end{align*}
by combining a double bracket $\Pi\colon A\otimes A\to A\otimes A$ which behaves like a Poisson bracket, a double derivation $\Theta\colon A\to A\otimes A$ which encodes the failure of $\Pi$ to be skew-symmetric, and the \textit{triple divergence map} which is a lift of a (double) divergence map associated with a connection $\nabla$ explained above. 
\end{comment}

More specifically, denoting the cyclic quotient of an algebra $A$ by $|A|$ and the space of \textit{double derivations} on $A$ by $\DDer(A)$, we will construct the triple divergence
\[
	\TDiv^\nabla\colon \DDer(A) \to |\mathrm{trip}(A)|\cong A\otimes |A| \oplus |A|\otimes A
\]
which is a lift of the (double) divergence map $\Div^\nabla$ associated with a connection $\nabla$ explained above. These maps are summarised in the following commutative diagram:
\[\begin{tikzcd}[cramped]
	\DDer(A) \arrow[r, "\TDiv^\nabla"] \arrow[d, "\mathrm{mult}_*"] & \text{\textbar}\trip(A)\text{\textbar} \arrow[d, "|\mathrm{mult}|"]\\
	\Der(A) \arrow[r, "\Div^\nabla"] & \text{\textbar}A\text{\textbar}\otimes \text{\textbar}A\text{\textbar}\,.
\end{tikzcd}\]
where vertical arrows are naturally induced from the multiplication map of $A$. By combining the triple divergence with a double bracket $\Pi\colon A\otimes A\to A\otimes A$, which is a map that behaves like a Poisson bracket, we obtain the following composition:
\begin{align*}
	\phi_{\Pi,\nabla}\colon A &\xrightarrow{\Pi} \DDer(A) \xrightarrow{\TDiv^\nabla} |\mathrm{trip}(A)|.\\
	a &\,\mapsto \Pi(a,\cdot)
\end{align*}

The map $\phi_{\Pi,\nabla}$ can be exploited to give an algebraic description of the loop operation $\mu^\mathsf{fr}$ for a suitable framing $\mathsf{fr}$, which is our main result in this paper, as we now explain. Let $\mathscr{G}$ be the fundamental groupoid of a surface with a non-empty boundary and a given set of base points on the boundary. Then, we can consider the linear groupoid $\mathbb{K}\mathscr{G}$ over a field $\mathbb{K}$, the cyclic quotient $|\mathbb{K}\mathscr{G}|$, and the space of non-commutative $1$-forms $\Omega^1\mathbb{K}\mathscr{G}$, which is a finitely generated projective $\mathbb{K}\mathscr{G}$-bimodule. We have a double bracket $\kappa\colon \mathbb{K}\mathscr{G}\otimes \mathbb{K}\mathscr{G}\to \mathbb{K}\mathscr{G}\otimes \mathbb{K}\mathscr{G}$ defined by smoothing the intersection of two paths. The map in question is of the form
\[
	\mu^\mathsf{fr}\colon \mathbb{K}\mathscr{G} \to |\mathbb{K}\mathscr{G}|\otimes \mathbb{K}\mathscr{G}\oplus \mathbb{K}\mathscr{G}\otimes |\mathbb{K}\mathscr{G}|,
\]
defined similarly to $\kappa$, but by smoothing \textit{self-intersections} of a path, which leads to increased algebraic complexity. The following is the main result of this paper:

\begin{theorem*}[Theorem \ref{thm:groupoidmu}]
Let $\mathcal{C}$ be a free-generating system of $\mathscr{G}$, $\mathsf{fr}$ be a framing of the surface normalised near the base points, and $\nabla\!_\mathcal{C}$ the connection on $\Omega^1\mathbb{K}\mathscr{G}$ defined by $\nabla\!_\mathcal{C}((dc)c^{-1}) = 0$ for all $c\in\mathcal{C}$. If $\mathcal{C}$ and $\mathsf{fr}$ are suitably chosen, we have $-\phi_{\kappa,\nabla\!_\mathcal{C}} = \mu^\mathsf{fr}$.
\end{theorem*}

\noindent In \cite{akkn}, the corresponding map is defined in the case of a surface with one base point as the unique lift of a given double divergence using an auxiliary variable; our construction is based on linear categories to deal with multiple base points and is more direct with the help of a connection.

The map $\phi_{\Pi,\nabla}$ is further used to define the \textit{modular vector field}: given a double derivation $\Theta\colon A\to A\otimes A$ encoding the failure of $\Pi$ to be skew-symmetric, we can define the modular vector field $\mathbf{m}_{\Pi,\nabla,\Theta}$ associated with the data $(\Pi,\nabla,\Theta)$ by the formula
\begin{align*}
	\mathbf{m}_{\Pi,\nabla,\Theta}\colon A&\to A\otimes |A|\\
	a& \mapsto \mathrm{fd}(\phi_{\Pi,\nabla}(a)) + \Theta'(a)\otimes |\Theta''(a)|
\end{align*}
where $\mathrm{fd}$ is a ``fold'' map. This turns out to be a derivation; see Definition-Lemma \ref{deflem:modvec} for the details. This can be used to give a sufficient condition for the skew-symmetry of the map $\delta^{\mathrm{Ham}_\Pi,\nabla}= \Div^\nabla\!\circ\,\mathrm{Ham}_\Pi$:

\begin{theorem*}[Theorem \ref{thm:modvec}]
Suppose that $\Theta = \ad_e$ for some symmetric $e$: $e'\otimes e'' = e''\otimes e' \in A\otimes A$. Then, the map $\delta^{\mathrm{Ham}_\Pi,\nabla}$ is skew-symmetric if and only if $\mathrm{sym}\circ|\mathbf{m}_{\Pi,\nabla,\ad_e}| = 0$. In particular, if $\mathbf{m}_{\Pi,\nabla,\ad_e} = 0$, $\delta^{\mathrm{Ham}_\Pi,\nabla}$ is skew-symmetric.
\end{theorem*}
\noindent This is a general case of $\delta^{\sigma,\nabla}$ discussed above since the map $\sigma$ is exactly equal to $\mathrm{Ham}_\kappa$, the Hamiltonian flow specified with the double bracket $\kappa$. More details are covered in Section \ref{sec:ncgeom}.\\

In the case of free associative algebras, double brackets in relation to the (linearised) Goldman Lie algebra are studied in the paper \cite{naef} by Naef with the flavour of geometry over the Lie operad. As an application, he obtained a Poisson morphism connecting the reduction of the space of coadjoint orbits of $\mathfrak{g} = \mathrm{Lie}(G)$ to the moduli space of flat $G$-bundles on a surface for a Lie group $G$, on which Goldman discovered a symplectic structure and hence the Goldman bracket: a loop operation based on intersections of loops.\\

\noindent\textbf{Organisation of the paper.} In Section \ref{sec:ncgeom}, we recall basic notions from non-commutative geometry. Section \ref{sec:div} is devoted to the construction of divergence maps and modular vector fields. Section \ref{sec:lincat} is a review of linear categories and modules over them. Finally, in Section \ref{sec:grpdmu}, we recall the definition of the loop operation $\mu^\mathsf{fr}$ and prove the theorem above.\\

\noindent\textbf{Acknowledgements.} The author is deeply thankful to Florian Naef for providing the key idea of considering modular vector fields in a non-commutative setting, to Nariya Kawazumi for carefully reading a draft of this paper, to Gw\'ena\"el Massuyeau for pointing out additional references, and to Geoffrey Powell and anonymous referees for improving many details.\\

\noindent\textbf{Conventions.} $\mathbb{K}$ is a field of characteristic zero. Unadorned tensor products are always over $\mathbb{K}$. We use Sweedler-type notation: an element of the tensor product $X\otimes X$ is denoted as $x = x'\otimes x''$.
\\

\section{Preliminaries on Non-commutative Geometry}\label{sec:ncgeom}

In this section, we recall some definitions from non-commutative geometry. See \cite{ginz} for the general theory of non-commutative geometry, and see \cite{doublep} for the theory of double brackets.

\begin{remark}
In \cite{doublep}, van den Bergh defined double brackets to be skew-symmetric. Our definition is based on \cite{akkn}, and the skew-symmetry is not imposed.
\end{remark}

Let $A$ be a unital associative $\mathbb{K}$-algebra, $S$ a $\mathbb{K}$-subalgebra of $A$, and $A^\mathrm{e} = A\otimes A^\mathrm{op}$ its enveloping algebra. We identify $A$-bimodules, left $A^\mathrm{e}$-modules and right $A^\mathrm{e}$-modules in the usual way. \\

\noindent\textbf{Notation.} For a unital associative $\mathbb{K}$-algebra $A$,  an element of the opposite algebra $A\op$ corresponding to $a\in A$ is denoted by $\bar a$ so that $\bar a\bar b = \overline{ba}$. The multiplication map of $A$ is denoted $\mathrm{mult}\colon A\otimes_SA\to A$.\\

First, $A\otimes A$ is endowed with the following two $A$-bimodule structures. The \textit{outer structure} is given by
	\begin{align*}
		x\cdot (a\otimes b)\cdot y = xa\otimes by,
	\end{align*}
	while the \textit{inner structure} is given by
	\begin{align*}
		x\cdot_\mathrm{in} (a\otimes b)\cdot_\mathrm{in} y = ay\otimes xb.
	\end{align*}
	Unless otherwise specified, we use the outer structure. With these structures, we have an isomorphism of $A^\mathrm{e}$-bimodules given by
	\begin{align*}
		(A\otimes A, \,\cdot\,, \,\cdot_\mathrm{in}) \cong (A^\mathrm{e},\textrm{left multiplication},\textrm{right multiplication}): a\otimes b \mapsto a\otimes\bar b\,.
	\end{align*}

\begin{definition}\ 
\begin{itemize}
	\item $\Omega^1_SA = \Ker(A\otimes_S A \xrightarrow{\mathrm{mult}} A)$ is the space of $1$\textit{-forms} relative to $S$, which is an $A$-bimodule. We put $da = 1\otimes a - a\otimes 1$. With this notation, we have the Leibniz rule for $d$:
	\[
		d(ab) = (da)b + a(db) \textrm{ for all } a,b \in A.
	\]	
	 We drop $S$ in the absolute case $S = \mathbb{K}$.
	\item  The algebra of \textit{differential forms} is defined by
	\[
		\Omega^\bullet_SA := T_A\Omega_S^1A = \bigoplus_{p\geq 0}\, \underbrace{\Omega^1_SA\otimes_A \cdots \otimes_A \Omega^1_SA}_{p\textrm{ times}}.
	\]
	Its typical element is of the form $a_0da_1\cdots da_p$ with the tensor symbols omitted, and we can always rewrite any element into such a form by the Leibniz rule. The \textit{exterior derivative} $d$ is the degree-$1$ differential (with Koszul signs) on $\Omega^\bullet_SA$ given by
	\[
		d(a_0da_1\cdots da_p) = da_0da_1\cdots da_p,
	\]
	which makes $\Omega^\bullet_SA$ into a differential graded algebra.
	\item For an $A$-bimodule $M$, a $\mathbb{K}$-linear map $f\colon A\to M$ is called an $S$\textit{-linear derivation} if it satisfies
	\begin{align*}
		f(ab)& = f(a)\cdot b + a\cdot f(b)\quad\textrm{ for } a,b\in A,\textrm{ and}\\
		f(s) &= 0 \hspace{78pt} \textrm{ for } s\in S\,.
	\end{align*}
	The space of all such maps is denoted by $\Der_S(A,M)$. We set $\Der_S(A) = \Der_S(A,A)$ and $\DDer_S(A) = \Der_S(A,A\otimes A)$. An element of the latter is called a \textit{double derivation}.
		
	\hspace{15pt}The functor $\Der_S(A,\,\cdot\,)$ is isomorphic to $\Hom_{A^\mathrm{e}}(\Omega^1_SA, \,\cdot\,)$, and the element corresponding to $f$ is denoted by $i_f$ with the relation $i_f\circ d = f$. In this sense, $d\colon A \to \Omega^1_SA$ is the universal $S$-linear derivation. The right $A^\mathrm{e}$-module structure on $\DDer_S(A)$ is induced from the inner structure on $A\otimes A$. On the other hand, $\Der_S(A)$ is not an $A$-(bi)module in general.
	\item For $f\in \Der_S(A)$, the \textit{contraction map} $i_f$ is the degree $(-1)$ derivation on $\Omega^\bullet_SA$ given by canonically extending the above $i_f$ defined on $1$-forms. The \textit{Lie derivative} $L_f$ is defined by Cartan's formula: $L_f = [d,i_f] = d\circ i_f + i_f\circ d$, which is a degree $0$ derivation. For example, we have
	\begin{align*}
		L_f(a_0da_1) = f(a_0)da_1 + a_0df(a_1).
	\end{align*}
	\item For an $A$-bimodule $M$, we put $|M| = M/\langle am-ma \mid a\in A, m\in M\rangle_{\mathbb{K}\textrm{-vect}}$. This is called the \textit{trace space} of $M$. If $f\colon A\to M$ is a ($\mathbb{K}$-linear) derivation, it descends to $|f|\colon |A| \to |M|$.\\
\end{itemize}
\end{definition}

\begin{definition}\ 
\begin{itemize}
	\item A left $A$-module $N$ is \textit{dualisable} if it is finitely generated and projective. This is equivalent to one of the following:
	\begin{enumerate}[(i)]
		\item For any left $A$-module $L$, the canonical map $N^*\otimes_{A}L\to \Hom_{A}(N,L)$ is an isomorphism. Here $N^*= \Hom_A(N,A)$ is the $A$-dual of $N$, which is a right $A$-module.
		\item For some non-negative integer $r$, we can take two sets of elements $e_1, \dotsc, e_r \in N$ and $e^*_1,\dotsc, e^*_r\in N^*$ such that $\sum_i e^*_i\cdot e_i = \id_N$.
	\end{enumerate}
	\item For an $A$-bimodule $M$ and a dualisable left $A$-module $N$, the trace map is defined by the composition
	\[
		\Tr\colon \Hom_A(N, M\otimes_A N) \cong N^*\otimes_A M \otimes_A N \xrightarrow{|\mathrm{ev}|} |M|,
	\]
	where $\mathrm{|ev|}(\theta\otimes m\otimes n) = |\theta(n)m|$. Using $e_i$ and $e^*_i$ above, the trace of $\varphi\colon N \to M\otimes_AN$ is given explicitly by
	\[
		\Tr(\varphi) = \sum_i (\id_M\otimes \,e^*_i)(\varphi(e_i))\,.
	\]
	In the special case of $M = A$, it is the usual trace of an endomorphism.\\
\end{itemize}
\end{definition}

\begin{lemma}\label{lem:tracefunc}
For $A$-bimodules $M_1$, $M_2$ and an $A$-bimodule homomorphism $g\colon M_1\to M_2$, we have a commutative square
\[\begin{tikzcd}[cramped]
	\Hom_A(N, M_1\otimes_A N) \arrow[r, "\Tr"] \arrow[d, "g_*"] & \text{\textbar}M_1\text{\textbar} \arrow[d, "|g|"]\\
	\Hom_A(N, M_2\otimes_A N) \arrow[r, "\Tr"] & \text{\textbar}M_2\text{\textbar}\,.
\end{tikzcd}\]
\end{lemma}
\noindent Proof. For $\theta\otimes m\otimes n \in N^*\otimes_A M\otimes_A N$, we have
\begin{align*}
	(|g|\circ |\mathrm{ev}|)(\theta\otimes m\otimes n) = |g(\theta(n)m)| = |\theta(n)g(m)| = |\mathrm{ev}|(\theta\otimes g(m)\otimes n)\,.
\end{align*}
This completes the proof.\qed\\

\begin{definition}\ 
\begin{itemize}
	\item An $S$-linear \textit{double bracket} on $A$ is a $\mathbb{K}$-linear map $\Pi\colon A\otimes A\to A\otimes A$ such that, for $a,b,c\in A$ and $s\in S$, 
	\begin{align*}
		\Pi(a,bc) &= \Pi(a,b)\cdot c + b\cdot \Pi(a,c)\,,\\
		\Pi(ab,c) &= \Pi(a,c)\cdot_\mathrm{in}b + a\cdot_\mathrm{in} \Pi(b,c)\,\textrm{ and}\\
		\Pi(a,s) &= \Pi(s,b) = 0\,.
	\end{align*}
	The first and the third conditions ensure that $\Pi(a,\cdot)$ is an $S$-linear double derivation on $A$, and the second and the third conditions say that the map $\Pi\colon A \to \DDer_S(A): a\mapsto \Pi(a,\cdot)$ is itself an $S$-linear derivation. Therefore, an $S$-linear double bracket $\Pi$ is equivalent to an element $i_\Pi \in \Hom_{A^\mathrm{e}}(\Omega_S^1A,\DDer_S(A))$.
	\item For an $S$-linear double bracket $\Pi$ on $A$, its \textit{inversion} $\Pi^\circ\colon A\otimes A\to A\otimes A$ is defined by
	\[
		\Pi^\circ(a,b) = \Pi''(b,a)\otimes \Pi'(b,a)\,.
	\]
	Here, we used the Sweedler notation $\Pi = \Pi'\otimes \Pi''$. This is again an $S$-linear double bracket. $\Pi$ is \textit{skew-symmetric} if $\Pi + \Pi^\circ = 0$ holds, and \textit{symmetric} if $\Pi = \Pi^\circ$ holds.
\end{itemize}
\end{definition}

There is a natural map $\DDer_S(A)\xrightarrow{\mathrm{mult}_*}\Der_S(A)$ induced by the multiplication map on $A$, and this descends to $|\!\DDer_S(A)|\to \Der_S(A)$. In addition, since $\Pi\colon A \to \DDer_S(A): a\mapsto \Pi(a,\cdot)$ is a derivation, it also descends to $|\Pi|\colon |A| \to |\!\DDer_S(A)|$. 
\begin{definition}
We define $\mathrm{Ham}_\Pi$ by the composition
\begin{align*}
	\mathrm{Ham}_\Pi\colon|A| \xrightarrow{|\Pi|} |\!\DDer_S(A)| \xrightarrow{\mathrm{mult}_*} \Der_S(A),
\end{align*}
which is thought of as a non-commutative geometry analogue of a Hamiltonian flow.
\end{definition}

These maps are summarised in the following commutative diagram:
\[\begin{tikzcd}[cramped]
	A\arrow[r, "\Pi"] \arrow[d, ->>] & \DDer_S(A) \arrow[d, "\mathrm{mult}_*"]\phantom{\,.}\\
	\text{\textbar}A\text{\textbar} \arrow[r, "\mathrm{Ham}_\Pi"] & \Der_S(A)\,.
\end{tikzcd}\]\vspace{3pt}

\begin{example}\label{ex:db}
The following are important examples of double brackets:
\begin{enumerate}[(1)]
	\item \label{ex:dbgr} Let $W$ be a $\mathbb{K}$-vector space equipped with a pairing $\pair{\cdot}{\cdot}\colon W\otimes W\to \mathbb{K}$. Then, we get the $\mathbb{K}$-linear double bracket on the tensor algebra $T(W)$, again denoted by $\pair{\cdot}{\cdot}$, by extending the pairing with the Leibniz rules above. Explicitly, we have, for $u_i, w_j\in W$,
	\[
		\pair{u_1\cdots u_r}{w_1\cdots w_s} = \sum_{\substack{1\leq i\leq r\\1\leq j\leq s}} \pair {u_i}{w_j} w_1\cdots w_{j-1}u_{i+1}\cdots u_r\otimes u_1\cdots u_{i-1}w_{j+1}\cdots w_s\,.
	\]
	If the pairing is skew-symmetric, so is the induced double bracket.
	\item Let $\Sigma$ be a connected oriented compact surface with non-empty boundary, $\ast$ a base point on $\partial \Sigma$, and $\pi = \pi_1(\Sigma,\ast)$ the fundamental group of $\Sigma$. Then, there is a $\mathbb{K}$-linear double bracket on the group algebra $\mathbb{K}\pi$, which is independently introduced as $\{\!\!\{ \cdot\,,\cdot \}\!\!\}^\eta$ in \cite{mastur}, and as $\kappa$ in \cite{kappa}, respectively. Adopting the notation in \cite{kappa}, it is given as follows. First, fix an embedded short positive arc $\nu\colon [0,1]\to\partial \Sigma$ with $\nu(1) = \ast$, and put $\nu(0) = \bullet$\,, by which $\pi_1(\Sigma,\bullet)$ is identified with $\pi_1(\Sigma,\ast)$. For generic loops $\alpha$ based at $\bullet$, and $\beta$ at $\ast$, respectively, we define
	\[
		\kappa(\alpha,\beta) = \sum_{p\in\alpha\cap\beta} \mathrm{sign}(\alpha,\beta;p)\beta_{\ast p}\alpha_{p\bullet}\nu\otimes \nu^{-1}\alpha_{\bullet p}\beta_{p\ast}\,,
	\] 
	where $\mathrm{sign}(\alpha,\beta;p)$ is the local intersection number with respect to the orientation of $\Sigma$, and $\beta_{\ast p}$ is the path along $\beta$ from $\ast$ to $p$ and so on. Compositions of paths are read from left to right. This is \textit{not} skew-symmetric, and the defect is given explicitly as
	\begin{align}\label{eq:kappaskew}
		(\kappa + \kappa^\circ)(\alpha,\beta) = \alpha\otimes \beta + \beta\otimes \alpha - 1\otimes \alpha\beta - \beta\alpha\otimes 1\,.
	\end{align}
	In \cite{mastur}, $\{\!\!\{ \cdot\,,\cdot \}\!\!\}^\eta$ was constructed from Turaev's homotopy intersection form $\eta$ via the formalism of Fox pairings on a Hopf algebra and was shown to satisfy the quasi-Poisson identity.\\
	\hspace{15pt}The right-hand side of \eqref{eq:kappaskew} is, in a sense, an ``inner'' double bracket. Namely, set $e = 1\otimes 1\in A\otimes A$ for any unital algebra $A$, and consider the inner double derivation on $A$ by $e$:
	\[
		\ad_e\colon A \to A\otimes A\colon x \mapsto [x,e] = x\otimes 1 - 1\otimes x\,.
	\]
	Then, the inner $A^\mathrm{e}$-module derivation $\ad_{\ad_e}\colon A\to \DDer_\mathbb{K}(A)$ is exactly the right-hand side of \eqref{eq:kappaskew}. To see this, we compute
	\begin{align*}
		\ad_{\ad_e}(\alpha) = \alpha\cdot \ad_e -\ad_e\cdot\alpha
	\end{align*}
	so that 
	\begin{align*}
		\ad_{\ad_e}(\alpha)(\beta)&= \alpha\cdot_\mathrm{in} (\beta\otimes 1 - 1\otimes \beta) -(\beta\otimes 1 - 1\otimes \beta)\cdot_\mathrm{in}\alpha \\
		&= \beta\otimes\alpha - 1\otimes \alpha\beta - \beta\alpha\otimes 1 + \alpha\otimes \beta.
	\end{align*}
	In other words, $\kappa$ is skew-symmetric up to an inner double bracket. The defect in the skew-symmetry of $\eta$ was originally analysed in \cite{mastur}, where it is denoted by $\rho_1$, the inner Fox pairing associated with $1\in A$.\\
	\hspace{15pt}The double bracket $\pair{\cdot}{\cdot}$ above is an algebraic counterpart of $\kappa$. In fact, if $\Sigma$ has only one boundary component, letting $W = H_1(\Sigma;\mathbb{K})$ and $\pair{\cdot}{\cdot}$ be the intersection form on $\Sigma$ recovers the associated graded of $\kappa$ under a suitable filtration (see Section 3.3 of \cite{akkn}). We will see the groupoid version of $\kappa$ in Section \ref{sec:grpdmu}, which gives an $S$-linear double bracket for some separable algebra $S$.

\end{enumerate}
\end{example}

\begin{remark}
Suppose that the pairing $\pair\cdot\cdot$ is skew-symmetric. Then, the above double brackets define Lie algebra structures on the spaces $|T(W)|$ and $|\mathbb{K}\pi|$, and each $\mathrm{Ham}_\Pi$ is a Lie algebra homomorphism. The case of $|T(W)|$ is known as the Necklace Lie algebra, while the case of $|\mathbb{K}\pi|$ is the Goldman Lie algebra of the surface $\Sigma$.
\end{remark}

Following \cite{akkn}, we make the following:
\begin{definition}
We define the map $\sigma\colon |\mathbb{K}\pi| \to \Der_\mathbb{K}(\mathbb{K}\pi)$ to be the Hamiltonian flow $\mathrm{Ham}_\kappa$ associated with the above $\kappa$.\\
\end{definition}

Next, we recall the definition of a connection in non-commutative geometry. Let $B$ be another $\mathbb{K}$-algebra, $R$ a subalgebra and $M$ a left $B$-module.
\begin{definition}
An $R$-linear connection on $M$ is a $\mathbb{K}$-linear map $\nabla\colon M\to \Omega^1_RB\otimes_BM$ such that
\[
	\nabla(bm) = db \otimes m + b\nabla(m)\;\textrm{ for } b\in B\textrm{ and } m\in m.
\]
A connection is \textit{flat} if the curvature $\nabla^2\colon M\to \Omega^2_RB\otimes_BM$ vanishes.
\end{definition}

\begin{example}
We recall two flat $\mathbb{K}$-linear connections defined in the previous paper \cite{toyo}.
\begin{itemize}
	\item For a $\mathbb{K}$-vector space $W$, the connection $\nabla\!_W$ is defined on $\Omega^1T(W)$ by $\nabla\!_W(dw) = 0$ for all $w\in W$.
	\item The connection $\nabla\!_\mathcal{C}$ is defined on $\Omega^1\mathbb{K}F_n$ by $\nabla\!_\mathcal{C}((dx_i)x_i^{-1}) = 0 $, where $F_n$ is the free group with the free generating system $\mathcal{C} = (x_1,\dotsc,x_n)$.\\
\end{itemize}
\end{example}

\section{Divergences and Non-commutative Modular Vector Fields}\label{sec:div}

\subsection{Double and Triple Divergences}
The main object in the previous paper \cite{toyo} is the \textit{double divergence map} associated with a connection. We shall now briefly recall these.

\begin{definition}
Let $A$ be a $\mathbb{K}$-algebra and $S$ a $\mathbb{K}$-subalgebra of $A$. $A$ is said to be \textit{formally smooth relative to} $S$ if the space of relative $1$-forms $\Omega^1_SA$ is $A^\mathrm{e}$-dualisable.
\end{definition}

\begin{example}
The group algebra $\mathbb{K}F_n$ and the free associative algebra $T(W)$ are both formally smooth relative to $\mathbb{K}$. For the proof, see Section 3 of \cite{toyo}.
\end{example}

Assume that $A$ is formally smooth relative to $S$ and let $\nabla\colon \Omega^1_SA\to \Omega^1_{S^\mathrm{e}}A^\mathrm{e}\otimes_{A^\mathrm{e}}\Omega^1_SA$ be an $S^\mathrm{e}$-linear connection. Note that the difference $(i_{f^\mathrm{e}}\otimes\id)\circ\nabla - L_f$ is $A^\mathrm{e}$-linear, so taking the trace makes sense.

\begin{definition}\label{def:conn}
	The \textit{double divergence map} on $\Der_S(A)$ associated with $\nabla$ is defined by 
	\[
		\Div^\nabla\colon\mathfrak \Der_S(A)\to |A^\mathrm{e}|\colon f\mapsto \Tr((i_{f^\mathrm{e}}\otimes\id)\circ\nabla - L_f)\,.
	\]
Here we denote $f^\mathrm{e} = f\otimes \id + \id \otimes f$, which is an $S^\mathrm{e}$-linear derivation on $A^\mathrm{e}$. We set
\[
	\delta^{\mathrm{Ham}_\Pi,\nabla} = \Div^\nabla\!\circ\,\mathrm{Ham}_\Pi \colon |A| \to |A^\mathrm{e}|\,.
\]
\end{definition}

\begin{remark}
For a more general construction involving a derivation action on a module, see Section 4 of \cite{toyo}. The definition above differs by a factor of $(-1)$ from the one in \cite{toyo}.
\end{remark}

In the case of the surface $\Sigma$, as in Example \ref{ex:db} (2), we have the following:

\begin{theorem}[\cite{akkn}, \cite{toyo}]\label{thm:cob}
Let $\Sigma$ be a connected oriented compact surface with non-empty boundary, $\mathcal{C} = (\alpha_i,\beta_i,\gamma_j)$ the free-generating system of $\pi = \pi_1(\Sigma,\ast_0)$ in Figure \ref{fig:gensys}, and $\mathsf{fr}$ a framing on $\Sigma$ with $\mathrm{rot}^\mathsf{fr}(|c|)=0$ for all $c\in\mathcal{C}$. Then, $-\delta^{\sigma,\nabla\!_\mathcal{C}}$ is equal to the framed version of the Turaev cobracket $\delta^\mathsf{fr}$.
\end{theorem}

The purpose of this paper is to prove an analogous statement to the above for Turaev's operation $\mu^\mathsf{fr}$, which is a lift of $\delta^\mathsf{fr}$ along the projection $A\twoheadrightarrow |A|$. With this in mind, we will construct a lift of the double divergence map, the \textit{triple divergence map}, from a connection. At the end of this section, we will also discuss the skew-symmetry of $\delta^{\mathrm{Ham}_\Pi,\nabla}$; this property is a non-commutative geometry analogue of the vanishing of the divergence.

\begin{definition}
Let $S$ be a $\mathbb{K}$-subalgebra of $A$. For an $S$-linear double derivation $\Theta\colon A\to A\otimes A$, the $\mathbb{K}$-linear map $\Theta\op\colon A\op\to A\op\otimes A\op$ is defined by
\[
	\Theta\op(\bar a) = \overline{\Theta''(a)}\otimes \overline{\Theta'(a)}\,.
\]
Here we denote $\Theta(a) = \Theta'(a)\otimes \Theta''(a)$.\vspace{3pt}
\end{definition}

\begin{lemma}
Under the identification between $A$-bimodules and $A\op$-bimodules, we have an isomorphism of $A$-bimodules $\DDer_S(A) \cong \DDer_{S\op}(A\op)$ given by $\Theta\mapsto \Theta\op$.
\end{lemma}
\noindent Proof. We first check that $\Theta\op$ is a double derivation. For $a,b\in A$, we have
\begin{align*}
	\Theta\op(\bar a\bar b) &= \Theta\op(\overline{ba})\\
	&= \overline{\Theta''(ba)}\otimes \overline{\Theta'(ba)}\\
	&= \overline{\Theta''(b) a} \otimes \overline{\Theta'(b)} + \overline{\Theta''(a)}\otimes \overline{b \Theta'(a)}\\
	&= \bar a\, \overline{\Theta''(b)} \otimes \overline{\Theta'(b)} + \overline{\Theta''(a)}\otimes \overline{\Theta'(a)}\, \bar b.
\end{align*}
This shows that $\Theta\op$ is a double derivation. $S$-linearity is clear from the construction. Next, we check that the map is an $A$-bimodule homomorphism. For $a,b,x\in A$, we have
\begin{align*}
	(a\cdot\Theta\op\cdot b)(x) &= (\bar b\cdot\Theta\op\cdot \bar a)(x)\\
	&= {\Theta\op}'(x)\,\bar a\otimes \bar b\,{\Theta\op}''(x)\\
	&=  \overline{\Theta''(x)} \,\bar a\otimes  \bar b\, \overline{\Theta'(x)}\\
	&=  \overline{a\Theta''(x)}\otimes \overline{\Theta'(x)b}\\
	&= (a\cdot\Theta\cdot b)\op (x).
\end{align*}
Lastly, the map is bijective since we have $(\Theta\op)\op = \Theta$. This completes the proof.\qed\vspace{3pt}

\begin{definition}
The \textit{triple} of $A$, $\trip(A)$, is defined by 
	\[
		\trip(A) = ((A\otimes A)\otimes A\op) \oplus (A\otimes (A\op \otimes A\op))\,
	\]
	with the \textit{outer} $A^\mathrm{e}$-bimodule structure given by
	\begin{align*}
		(a\otimes\bar b)\cdot (x\otimes y\otimes \bar z) = ax\otimes y \otimes \bar b\bar z,\hspace{23pt} (x\otimes y\otimes \bar z)\cdot (a\otimes\bar b) = x\otimes ya\otimes \bar z\bar b\,,\\
		(a\otimes\bar b)\cdot (x\otimes \bar y\otimes \bar z) = ax\otimes \bar b\bar y \otimes \bar z,\textrm{ and }  (x\otimes \bar y\otimes \bar z)\cdot (a\otimes\bar b) = xa\otimes \bar y\otimes \bar z\bar b\,.
	\end{align*}
	This module structure is briefly denoted by dots:
	\[
		\trip(A) = ((.A\otimes A.)\otimes .A\op .) \oplus (.A.\otimes (.A\op \otimes A\op .))\,,
	\]
	indicating the places where multiplications in $A$ and $A\op$ are applied. In addition, the \textit{inner} structure on $\trip(A)$ is the $A$-bimodule structure given by
	\begin{align*}
		((A.\otimes .A)\otimes A\op)\oplus (A\otimes (A\op. \otimes .A\op))\,.\\
	\end{align*}
\end{definition}

\begin{lemma}\label{lem:tripM}
For a left $A^\mathrm{e}$-module $M$, we have a canonical isomorphism of left $A^\mathrm{e}$-modules:
	\begin{align*}
		\begin{aligned}
			\trip(A)\otimes_{A^\mathrm{e}}M &\cong (.A\otimes M.) \oplus (.M\otimes A.)\\
			(x\otimes y\otimes \bar z)\otimes m &\leftrightarrow \;x\otimes ymz\,,\\
			(x\otimes \bar y\otimes \bar z)\otimes m &\leftrightarrow \hspace{55pt} xmz\otimes  y\,.
		\end{aligned}
	\end{align*}
	Here, $\trip(A)$ is equipped with the outer structure.
\end{lemma}
\noindent Proof. This can be checked by direct verification.\qed\\

Recall that $\mathrm{mult}\colon A\otimes A \to A$ is an $A$-bimodule morphism with the outer structure on $A\otimes A$.
	
\begin{definition}
The multiplication map $\mathrm{mult}\colon \trip(A) \to A^\mathrm{e}$ is defined by
\begin{align*}
	&\mathrm{mult}(x\otimes y\otimes \bar z) = xy\otimes \bar z\,\textrm{ and}\\
	&\mathrm{mult}(x\otimes \bar y\otimes \bar z) = x\otimes \bar y\bar z\,.
\end{align*}
This is an $A^\mathrm{e}$-bimodule homomorphism.
\end{definition}

\begin{lemma}
We have a canonical isomorphism as $A$-bimodules:
\begin{align}\label{eq:tripA}\begin{aligned}
	|\trip(A)| &\cong (.A.\otimes |A|) \oplus (|A|\otimes .A .)\\
	x\otimes y\otimes \bar z &\leftrightarrow \;yx\otimes |z|\,,\\
	x\otimes \bar y\otimes \bar z &\leftrightarrow\hspace{61pt} |x|\otimes yz\,.
\end{aligned}\end{align}
Here, $|\trip(A)|$ is the trace space as an $A^\mathrm{e}$-bimodule, which is an $A$-bimodule by the inner structure. With this identification, the map $|\mathrm{mult}|\colon |\mathrm{trip}(A)| \to |A^\mathrm{e}|$ corresponds to
\begin{align*}
	|\cdot|\otimes\id + \id\otimes|\cdot|\colon A\otimes|A|\oplus |A|\otimes A&\to |A^\mathrm{e}|\\
	x\otimes |y| \hspace{47pt}&\mapsto |x|\otimes |\bar y|\\
	|x|\otimes y \hspace{2pt}&\mapsto |x|\otimes |\bar y|\,.
\end{align*}
\end{lemma}
\noindent Proof. This can be checked by direct verification.\qed\\

\begin{definition} Let $S$ be a $\mathbb{K}$-subalgebra of $A$ and $\Theta\in\DDer_S(A)$.
\begin{itemize}
	\item $\Theta^\mathrm{e}\colon A^\mathrm{e} \to \trip(A)$ is defined by $\Theta^\mathrm{e}(a\otimes \bar b) = \Theta(a)\otimes \bar b +  a\otimes \Theta\op(\bar b)$. This is an $S^\mathrm{e}$-linear $A^\mathrm{e}$-bimodule derivation with the structure on $\trip(A)$ given above. Therefore it induces the map $i_{\Theta^\mathrm{e}}\in\Hom_{(A^\mathrm{e})^\mathrm{e}}(\Omega_{S^\mathrm{e}}^1A^\mathrm{e},\trip(A))$.
	\item Let $M$ be a left $A^\mathrm{e}$-module and $\nabla\colon M\to \Omega_{S^\mathrm{e}}^1A^\mathrm{e}\otimes_{A^\mathrm{e}}M$ an $S^\mathrm{e}$-linear connection. Then, we define the map $i_{\Theta^\mathrm{e}}\nabla$ to be the following composition:
	\[
		i_{\Theta^\mathrm{e}}\nabla\colon M\xrightarrow{\nabla} \Omega_{S^\mathrm{e}}^1A^\mathrm{e}\otimes_{A^\mathrm{e}}M\xrightarrow{i_{\Theta^\mathrm{e}}\otimes\id} \trip(A)\otimes_{A^\mathrm{e}}M.
	\]
	\item For an $S$-linear double derivation $\Theta$, the map $L_{\Theta}\colon \Omega_S^1A \to (.\Omega_S^1A\otimes A.)\oplus (.A\otimes \Omega_S^1A.)$ is given by
	\[
		L_\Theta(a_0da_1) = \Theta'(a_0)\otimes \Theta''(a_0)da_1 + a_0d\Theta'(a_1)\otimes \Theta''(a_1) + a_0\Theta'(a_1)\otimes d\Theta''(a_1)\,.
	\]
	The target is canonically isomorphic to $\trip(A)\otimes_{A^\mathrm{e}} \Omega_S^1A$ by Lemma \ref{lem:tripM}.
\end{itemize}
\end{definition}

The map $L_\Theta$ is the double version of the Lie derivative: putting $f = \mathrm{mult}\circ\Theta \in \Der_S(A)$, we have  the following commutative diagram: 
\[\begin{tikzcd}[cramped]
	\Omega_S^1A \arrow[r, "L_\Theta"] \arrow[d, equal] & (\Omega_S^1A\otimes A)\oplus (A\otimes \Omega_S^1A) \arrow[d]\\
	\Omega_S^1A \arrow[r, "L_f"] & \Omega_S^1A
\end{tikzcd}\]
where the right arrow is the composition
\begin{align*}
	(\Omega_S^1A\otimes A)\oplus (A\otimes \Omega_S^1A) \twoheadrightarrow (\Omega_S^1A\otimes_A A)\oplus (A\otimes_A \Omega_S^1A) \cong \Omega_S^1A \oplus \Omega_S^1A \xrightarrow{\mathrm{add}} \Omega_S^1A.\\
\end{align*}
	
Now, we introduce the key ingredient in this paper:
\begin{definition}
Suppose that $A$ is formally smooth relative to $S$. We define the \textit{triple divergence} associated with an $S^\mathrm{e}$-linear connection $\nabla$ on $\Omega^1_SA$ by
\[
	\TDiv^\nabla\colon \DDer_S(A)\to|\trip(A)|\colon \Theta\to \Tr(i_{\Theta^\mathrm{e}}\nabla - L_{\Theta}),
\]
which is well-defined since the difference $i_{\Theta^\mathrm{e}}\nabla - L_{\Theta}$ is $A^\mathrm{e}$-linear.
\end{definition}

\begin{lemma}\label{lem:tdivlift}
We have the following commutative diagram:
\[\begin{tikzcd}[cramped]
	\DDer_S(A) \arrow[r, "\TDiv^\nabla"] \arrow[d, "\mathrm{mult}_*"] & \text{\textbar}\trip(A)\text{\textbar} \arrow[d, "|\mathrm{mult}|"]\\
	\Der_S(A) \arrow[r, "\Div^\nabla"] & \text{\textbar}A^\mathrm{e}\text{\textbar}\,.
\end{tikzcd}\]
Thus, $\TDiv^\nabla$ is a lift of $\Div^\nabla$ in this sense.
\end{lemma}
\noindent Proof. For $\Theta\in\DDer_S(A)$, put $f = \mathrm{mult}\circ \Theta \in\Der_S(A)$. It is sufficient to show that
\[
	(\mathrm{mult}\otimes \id)\circ(i_{\Theta^\mathrm{e}}\nabla - L_{\Theta}) = i_{f^\mathrm{e}}\nabla - L_f
\]
since we can apply Lemma \ref{lem:tracefunc} with $g = \mathrm{mult}$ to compute
\begin{align*}
	|\mathrm{mult}|(\TDiv^\nabla(\Theta)) &= (|\mathrm{mult}|\circ\Tr)(i_{\Theta^\mathrm{e}}\nabla - L_{\Theta})\\
	&= \Tr((\mathrm{mult}\otimes \id)\circ(i_{\Theta^\mathrm{e}}\nabla - L_{\Theta}))\\
	&= \Tr(i_{f^\mathrm{e}}\nabla - L_f)\\
	&= \Div^\nabla(f)\,.
\end{align*}
For the contraction $i_{\Theta^\mathrm{e}}$, we have, for $x,y,z,w\in A$,
\begin{align*}
	(\mathrm{mult}\circ i_{\Theta^\mathrm{e}})((x\otimes \bar y) d(z\otimes\bar w)) &= \mathrm{mult}((x\otimes \bar y) \Theta^\mathrm{e}(z\otimes\bar w))\\
	&= \mathrm{mult}((x\otimes \bar y) (\Theta(z)\otimes \bar w + z\otimes \Theta\op(\bar w)))\\
	&= (x\otimes \bar y) \cdot \mathrm{mult} (\Theta(z)\otimes \bar w + z\otimes \Theta\op(\bar w))\\
	&= (x\otimes \bar y) \cdot (f(z)\otimes \bar w + z\otimes \overline{f(w)})\\
	&= i_{f^\mathrm{e}}((x\otimes \bar y) d(z\otimes\bar w))\,.
\end{align*}
For the Lie derivative $L_\Theta$, we have, for $x,y\in A$,
\begin{align*}
	(\mathrm{mult}\circ L_\Theta)(xdy) &= \mathrm{mult}( \Theta'(x)\otimes \Theta''(x)dy + xd\Theta'(y)\otimes \Theta''(y) + x\Theta'(y)\otimes d\Theta''(y))\\
	&= \Theta'(x)\Theta''(x)dy + xd\Theta'(y) \Theta''(y) + x\Theta'(y) d\Theta''(y)\\
	&= f(x)dy + xdf(y)\\
	&= L_f(xdy).
\end{align*}
This completes the proof. \qed\\

\begin{proposition}\label{prop:multprop}
$\TDiv^\nabla$ has the following multiplicative property: for $\Theta\in\DDer_S(A)$ and $a,b\in A$, 
\begin{align*}
	\TDiv^\nabla((a\otimes\bar b)\cdot \Theta) = (a\otimes \bar b)\cdot \TDiv^\nabla(\Theta) - |\Theta'(b)|\otimes a\Theta''(b) - \Theta'(a)b\otimes |\Theta''(a)|\,.
\end{align*}
\end{proposition}
\noindent Proof. First of all, we have, in $(\Omega_S^1A\otimes A)\oplus (A\otimes \Omega_S^1A)$,
\[
	L_{(a\otimes\bar b)\cdot  \Theta}(\omega) =  L'_\Theta(\omega)b\otimes  aL''_\Theta(\omega) + i'_\Theta(\omega)b\otimes da\, i''_\Theta(\omega) + i'_\Theta(\omega)db\otimes ai''_\Theta(\omega)
\]
for $a,b\in A$ and $\omega\in\Omega_S^1A$. In addition, we have 
\[
	i_{((a\otimes\bar b)\cdot \Theta)^\mathrm{e}}\nabla = (a\otimes\bar b)\cdot_{\mathrm{in}}  i_{\Theta^\mathrm{e}}\nabla
\]
where $\cdot_\mathrm{in}$ indicates the action by the inner structure.
Since $M = \Omega_S^1A$ is dualisable by the assumption, we take a set of elements $e_1,\dotsc,e_r\in M$ and $e^*_1,\dotsc,e^*_r\in M^*$ with $\id_M = \sum_j e_j^*\cdot e_j$. Then, we compute the trace:
\begin{align*}
	\TDiv^\nabla((a\otimes\bar b)\cdot \Theta) &= \Tr(i_{((a\otimes\bar b)\cdot \Theta)^\mathrm{e}}\nabla - L_{(a\otimes\bar b)\cdot \Theta})\\
	&= \Big|\sum_j (\id_{\trip(A)}\otimes \,e^*_j)(i_{((a\otimes\bar b)\cdot \Theta)^\mathrm{e}}\nabla(e_j) - L_{(a\otimes\bar b)\cdot \Theta}(e_j))\Big|\\
	&= \Big|\sum_j (\id_{\trip(A)}\otimes \,e^*_j)\Big((a\otimes\bar b)\cdot_{\mathrm{in}}  i_{\Theta^\mathrm{e}}\nabla(e_j) \\
	&\qquad \qquad - (L'_\Theta(e_j)b\otimes  aL''_\Theta(e_j) + i'_\Theta(e_j)b\otimes da\, i''_\Theta(e_j) + i'_\Theta(e_j)db\otimes ai''_\Theta(e_j))\Big)\Big|\\
	&= (a\otimes\bar b)\cdot\Big|\sum_j (\id_{\trip(A)}\otimes \,e^*_j)(  i_{\Theta^\mathrm{e}}\nabla(e_j) - L_\Theta(e_j))\Big|\\
	&\qquad - \sum_j \Big( {e_j^*}'(da) i'_\Theta(e_j)b\otimes |{e^*_j}''(da)i''_\Theta(e_j)| + |i'_\Theta(e_j){e_j^*}''(db)|\otimes ai''_\Theta(e_j){e_j^*}''(db)\Big)\,.
\end{align*}
Applying $i_\Theta$ to the equality $da = \sum_j e_j^*(da)e_j$, we have
\begin{align*}
	\Theta(a) = \sum_j e_j^*(da)\cdot i_\Theta(e_j) = \sum_j {e_j^*}'(da)i'_\Theta(e_j) \otimes i''_\Theta(e_j){e_j^*}''(da)\,.
\end{align*}
Hence we obtain
\[
	\TDiv^\nabla((a\otimes\bar b)\cdot \Theta) = (a\otimes \bar b)\cdot \TDiv^\nabla(\Theta) - |\Theta'(b)|\otimes a\Theta''(b) - \Theta'(a)b\otimes |\Theta''(a)|\,,
\]
which completes the proof.\qed

\subsection{Modular Vector Fields}
Now, we introduce a modular vector field, whose purpose is to measure to what extent $\Pi$ and $\nabla$ are compatible. The definition is given in Definition-Lemma \ref{deflem:modvec}. It is an analogue of a modular vector field in Poisson geometry, so we quickly recall the construction on a manifold. Let $X$ be a smooth oriented Poisson manifold and $\mu$ a volume form. The \textit{modular vector field} associated with $\mu$ is the composition
\[
	\mathbf{m}_\mu\colon C^\infty(X)\xrightarrow{\mathrm{Ham}} \Der_\mathbb{R}(C^\infty(X)) \xrightarrow{\sdiv_\mu} C^\infty(X)\,,
\]
which is a vector field by the skew-symmetry of the Poisson bivector. The manifold $X$ is called \textit{unimodular} if $\mathbf{m}_\mu$ vanishes for some $\mu$, and this holds when, for example, the Poisson structure is induced from a symplectic structure; this is highly relevant to the case of surfaces considered above since the Goldman Lie algebra is a non-commutative counterpart of the symplectic structure on the moduli space of flat principal bundles.

With this seen, we define the following map, but using the lift $\TDiv^\nabla$:

\begin{definition}\label{def:phipinabla}
Let $\Pi\colon A\to\DDer_S(A)$ be an $S$-linear double bracket on $A$ and $\nabla$ as above. We put
\[
	\phi_{\Pi,\nabla}\colon A\xrightarrow{\Pi}\DDer_S(A)\xrightarrow{\TDiv^\nabla} |\trip(A)|\,.
\]
\end{definition}

\begin{lemma}\label{lem:multprop}
For $a,b\in A$, we have
\[
	\phi_{\Pi,\nabla}(ab) = \phi_{\Pi,\nabla}(a)\cdot b + a\cdot \phi_{\Pi,\nabla}(b) - (|\cdot|\otimes \id)\Pi(a,b) - (\id\otimes|\cdot|)\Pi(b,a)\,.
\]
Therefore $\phi_{\Pi,\nabla}$ induces $|\phi_{\Pi,\nabla}|\colon |A| \to |A^\mathrm{e}|$ which makes the following diagram commute:
\[\begin{tikzcd}[cramped]
	A\arrow[r, "\phi_{\Pi,\nabla}"] \arrow[d, ->>] & \text{\textbar}\trip(A)\text{\textbar} \arrow[d, "|\mathrm{mult}|"]\\
	\text{\textbar}A\text{\textbar} \arrow[r, "|\phi_{\Pi,\nabla}|"] & \text{\textbar}A^\mathrm{e}\text{\textbar}\,.
\end{tikzcd}\]
\end{lemma}
\noindent Proof. This follows from Proposition \ref{prop:multprop}.\qed

\begin{lemma}
The map $|\phi_{\Pi,\nabla}|$ is equal to $\delta^{\mathrm{Ham}_\Pi,\nabla}$ defined in Definition \ref{def:conn}.
\end{lemma}
\noindent Proof. Combining the definition of $\mathrm{Ham}_\Pi$ with Lemma \ref{lem:tdivlift}, we have the following commutative diagram:
\[\begin{tikzcd}[cramped]
	A\arrow[r, "\Pi"] \arrow[d, ->>] & \DDer_S(A) \arrow[r, "\TDiv^\nabla"] \arrow[d, "\mathrm{mult}_*"] & \text{\textbar}\trip(A)\text{\textbar} \arrow[d, "|\mathrm{mult}|"]\\
	\text{\textbar}A\text{\textbar} \arrow[r, "\mathrm{Ham}_\Pi"] & \Der_S(A) \arrow[r, "\Div^\nabla"] & \text{\textbar}A^\mathrm{e}\text{\textbar}\,.
\end{tikzcd}\]
Since the first row is $\phi_{\Pi,\nabla}$ and the second row is $\delta^{\mathrm{Ham}_\Pi,\nabla}$, we have $\delta^{\mathrm{Ham}_\Pi,\nabla} = |\phi_{\Pi,\nabla}|$ as desired.\qed\\

By Lemma \ref{lem:multprop}, $\phi_{\Pi,\nabla}$ is generally not a derivation. However, a slight modification will do the job if $\Pi$ is nearly skew-symmetric. This corresponds to the situation where we needed the skew-symmetry of the Poisson bivector on a manifold $X$. We first define the fold map by
\begin{align*}
	\mathrm{fd}\colon |\trip(A)|\cong  (A\otimes |A|) \oplus (|A|\otimes A) &\to A\otimes |A|\,,\\
	 x\otimes |y| + |z|\otimes w \hspace{10.5pt}&\mapsto x\otimes |y| + w\otimes |z|\,.
\end{align*}

\begin{deflem}\label{deflem:modvec}
Let $\Pi$ be a double bracket on $A$ with $\Pi + \Pi^\circ = \ad_\Theta$ for some $\Theta\in \DDer_S(A)$. Then, the map\\[-25pt]
\begin{align*}
	\mathbf{m}_{\Pi,\nabla,\Theta}\colon A&\to A\otimes |A|\\
	a& \mapsto \mathrm{fd}(\phi_{\Pi,\nabla}(a)) + \Theta'(a)\otimes |\Theta''(a)|
\end{align*}
is an $S$-linear $A$-bimodule derivation. We call it the modular vector field associated with the triple $(\Pi,\nabla,\Theta)$.
\end{deflem}
\noindent Proof. For $a,b\in A$, we have
\begin{align*}
	\mathbf{m}_{\Pi,\nabla,\Theta}(ab) &=  \mathrm{fd}(\phi_{\Pi,\nabla}(ab)) + \Theta'(ab)\otimes |\Theta''(ab)|\\
	&= \mathrm{fd}\Big(\phi_{\Pi,\nabla}(a)\cdot b + \phi_{\Pi,\nabla}(b) - (|\cdot|\otimes \id)\Pi(a,b) - (\id\otimes|\cdot|)\Pi(b,a)\Big) + \Theta'(ab)\otimes |\Theta''(ab)|\\
	&= \mathrm{fd}\Big(\phi_{\Pi,\nabla}(a)\cdot b + a\cdot\phi_{\Pi,\nabla}(b)\Big) - (\id\otimes|\cdot|)(\Pi+\Pi^\circ)(b,a) + \Theta'(ab)\otimes |\Theta''(ab)|\\
	&= \mathrm{fd}\Big(\phi_{\Pi,\nabla}(a)\cdot b + a\cdot\phi_{\Pi,\nabla}(b)\Big) \\
	&\qquad - ( \Theta'(a)\otimes |b\Theta''(a)| - \Theta'(a)b\otimes |\Theta''(a)|) + \Theta'(a)\otimes |\Theta''(a)b| + a\Theta'(b)\otimes |\Theta''(b)|\\
	&=\mathrm{fd}\Big(\phi_{\Pi,\nabla}(a)\cdot b + a\cdot\phi_{\Pi,\nabla}(b)\Big) + \Theta'(a)b\otimes |\Theta''(a)| + a\Theta'(b)\otimes |\Theta''(b)|\\
	&= \mathbf{m}_{\Pi,\nabla,\Theta}(a) \cdot b + a\cdot \mathbf{m}_{\Pi,\nabla,\Theta}(b) \,.
\end{align*}
Furthermore, $\mathbf{m}_{\Pi,\nabla,\Theta}$ is $S$-linear since $\Pi$ and $\Theta$ are. This completes the proof. \qed\\

Finally, we discuss the skew-symmetry of the map $\delta^{\mathrm{Ham}_\Pi,\nabla}$ under a mild hypothesis.

\begin{theorem}\label{thm:modvec}
Let $\Pi$ be a double bracket on $A$ with $\Pi + \Pi^\circ = \ad_{\ad_e}$ for some symmetric $e$: $e'\otimes e'' = e''\otimes e'$. Then, the map $\delta^{\mathrm{Ham}_\Pi,\nabla}$ is skew-symmetric if and only if $\mathrm{sym}\circ|\mathbf{m}_{\Pi,\nabla,\ad_e}| = 0$. In particular, if $\mathbf{m}_{\Pi,\nabla,\ad_e} = 0$, $\delta^{\mathrm{Ham}_\Pi,\nabla}$ is skew-symmetric.
\end{theorem}
\noindent Proof. First, $\delta^{\mathrm{Ham}_\Pi,\nabla}$ is skew-symmetric if and only if the symmetrised version $\mathrm{sym}\circ|\phi_{\Pi,\nabla}|$ vanishes. On the other hand, we have
\begin{align*}
	\mathrm{sym}\circ|\mathbf{m}_{\Pi,\nabla,\ad_e}| &= \mathrm{sym}\circ|\mathrm{fd}\circ \phi_{\Pi,\nabla}|+ \mathrm{sym}(|\ad_e(a)|)\,,
\end{align*}
and each term is further computed as follows:
\begin{align*}
	&\mathrm{sym}\circ |\mathrm{fd}\circ \phi_{\Pi,\nabla}| = \mathrm{sym}\circ |\phi_{\Pi,\nabla}|\textrm{ and}\\
	&\mathrm{sym}(|\ad_e(a)|) = \mathrm{sym}(|ae'|\otimes |e''| - |e'|\otimes |ae''|) = \mathrm{sym}(|ae'|\otimes |e''| -|e''|\otimes |ae'|) = 0\,.
\end{align*}
Therefore, we have $\mathrm{sym}\circ |\phi_{\Pi,\nabla}| = \mathrm{sym}\circ|\mathbf{m}_{\Pi,\nabla,\ad_e}|$. This completes the proof.\qed\\

\begin{remark}
The assumption that $e$ is symmetric is reasonable since $\Pi+\Pi^\circ$ is. Also, it is enough to check the condition $\mathbf{m}_{\Pi,\nabla,\ad_e} = 0$ on generators of $A$ as an algebra since it is a derivation by Definition-Lemma \ref{deflem:modvec}.
\end{remark}

\begin{example}
In the setting of Example \ref{ex:db} (1) with a skew-symmetric pairing, we have, for $u,w\in W$,
\[
	L_{\mathrm{Ham}_{\pair u\cdot}}(dw) = \pair uw\cdot d(1\otimes 1) = 0\;\textrm{ and }\; \nabla\!_W(dw) = 0\,.
\]
Hence the map $\phi_{\pair{\cdot}{\cdot},\nabla\!_W}(u)$ vanishes for all $u\in W$. Since we can take $e=0$ in this case, we have $\mathbf{m}_{\pair{\cdot}{\cdot},\nabla\!_W,\ad_e}=0$ and $\delta^{\mathrm{Ham}_{\pair{\cdot}{\cdot}},\nabla\!_W}$ is skew-symmetric.
\end{example}
\noindent We will discuss the case of the double bracket $\kappa$ in the more general setting of groupoids in the later sections.

\begin{remark}
The map $\TDiv$ in the case of $A = T(W)$ is the same as the maps $\mathfrak{b}_r+\mathfrak{b}_l$ in Theorem 4.50 of \cite{akkn}. They defined it as the unique extension of the divergence map to the space of double derivations by introducing an auxiliary variable $s$, whereas we constructed it directly using a connection.
\end{remark}

Differential-operator-like maps (such as derivations and divergences here) in a non-commutative setting are systematically studied in the paper \cite{ncdiffop} by Ginzburg and Schedler, as operators on the ``Fock space'' $\mathcal{F}(A)$ associated with $A$; compare Equation (1.4.4) there with Lemma \ref{lem:multprop} here.\\

\section{Non-commutative Geometry on a Linear Category} \label{sec:lincat}

In this section, we consider non-commutative geometry on linear categories, which is eventually reduced to the relative version on algebras discussed in the sections above. Another similar approach using a groupoid can be seen in the paper \cite{grpd} by Kawazumi and Kuno.

We start with basic notions to fix the notation.

\begin{definition}\ 
\begin{itemize}
	\item $\mathbf{Vect}_\mathbb{K}$ is the category of $\mathbb{K}$-vector spaces and $\mathbb{K}$-linear maps.
	\item A \textit{$\mathbb{K}$-linear category} is a $\mathbf{Vect}_\mathbb{K}$-enriched small category. Namely, a $\mathbb{K}$-linear category $\mathscr{A}$ consists of:
	\begin{enumerate}[(i)]
		\item the set $\Ob(\mathscr{A})$ of objects,
		\item the $\mathbb{K}$-vector space $\mathscr{A}(v,w)$ of morphisms from $v$ to $w$ for each pair of objects $(v,w)$, and
		\item the identity morphism $1_v\in \mathscr{A}(v,v)$ for each $v\in\mathrm{Ob}(\mathscr{A})$;
	\end{enumerate}
	with $\mathbb{K}$-linear composition map
	\begin{align*}
		\mathscr{A}(u,v)\otimes \mathscr{A}(v,w) &\to \mathscr{A}(u,w):\\
		(u\xrightarrow{x}v)\otimes (v\xrightarrow{y}w)& \mapsto (u\xrightarrow{xy}w).
	\end{align*}
	Here, the composition should be read from left to right so that it is compatible with the multiplication in the fundamental groupoid of a surface. Morphisms between $\mathbb{K}$-linear categories are $\mathbb{K}$-linear functors. The category of $\mathbb{K}$-linear categories is denoted by $\mathbf{Cat}(\mathbf{Vect}_\mathbb{K})$.
	\item The tensor product of $\mathbb{K}$-linear categories is taken in $\mathbf{Cat}(\mathbf{Vect}_\mathbb{K})$. Namely, for $\mathbb{K}$-linear categories $\mathscr{A}$ and $\mathscr{A}'$, $\mathscr{A}\otimes\mathscr{A}'$ is defined by $\Ob(\mathscr{A}\otimes\mathscr{A}') = \Ob(\mathscr{A})\times \Ob(\mathscr{A}')$ and $(\mathscr{A}\otimes\mathscr{A}')((v,v'),(w,w')) = \mathscr{A}(v,w)\otimes \mathscr{A}'(v',w')$ with the natural composition map.
	\item A \textit{left $\mathscr{A}$-module} $\mathscr{N}$ is a $\mathbb{K}$-linear functor $\mathscr{N}\colon\mathscr{A}\op \to \mathbf{Vect}_\mathbb{K}$. Namely, it is a collection of $\mathbb{K}$-vector spaces $\mathscr{N}(v)$ indexed by $v\in\Ob(\mathscr{A})$, and the action
	\[
		\mathscr{A}(w,v)\otimes \mathscr{N}(v) \to \mathscr{N}(w)\colon a\otimes n \to an
	\]
	satisfying the associativity and unit conditions. A morphism between left $\mathscr{A}$-modules is a natural transformation. A right $\mathscr{A}$-module is a functor $\mathscr{A} \to \mathbf{Vect}_\mathbb{K}$, and an $\mathscr{A}$-bimodule is a functor $\mathscr{A}^\mathrm{e}=\mathscr{A}\otimes\mathscr{A}\op \to \mathbf{Vect}_\mathbb{K}$.\\
	\hspace{15pt}For the sake of clarity, we write explicitly the definition of $\mathscr{A}$-bimodule $\mathscr{M}$: it is a collection of $\mathbb{K}$-vector spaces $\mathscr{M}(v,w)$ indexed by $v,w\in\Ob(\mathscr{A})$, and the action
	\[
		\mathscr{A}(v_1,v_2)\otimes \mathscr{M}(v_2,v_3)\otimes \mathscr{A}(v_3,v_4) \to \mathscr{M}(v_1,v_4)\colon x\otimes m \otimes y \to xmy
	\]
	satisfying the associativity and unit conditions, such that the left and right actions commute.
	\item Let $\mathscr{M}$ be an $\mathscr{A}$-bimodule. A \textit{derivation} $f\colon\mathscr{A}\to\mathscr{M}$ is a collection of $\mathbb{K}$-linear maps
	\[
		f(v,w)\colon \mathscr{A}(v,w)\to \mathscr{M}(v,w)
	\]
	for $v,w\in\Ob(\mathscr{A})$ satisfying the Leibniz rule $f(ab) = f(a)b + af(b)$ for any composable pair $(a,b)$. Denote by $\Der(\mathscr{A},\mathscr{M})$ the space of derivations.
	\item A \textit{double derivation} on $\mathscr{A}$ is a derivation of the form $\mathscr{A}\to \mathscr{A}\otimes \mathscr{A}$.
	\item A double bracket $\Pi\colon \mathscr{A}\otimes\mathscr{A}\to \mathscr{A}\otimes \mathscr{A}$ is a collection of $\mathbb{K}$-linear maps
	\[
		\Pi(v_1,v_2,v_3,v_4)\colon \mathscr{A}(v_1,v_2)\otimes \mathscr{A}(v_3,v_4) \to\mathscr{A}(v_3,v_2)\otimes \mathscr{A}(v_1,v_4)
	\]
	satisfying the two Leibniz rules as before for each composable set of maps.
\end{itemize}
\end{definition}

Below is a principal example of a projective module.

\begin{proposition}\label{prop:principaldualisable}
Let $\mathscr{A}$ be a $\mathbb{K}$-linear category. For $v\in\Ob(\mathscr{A})$, the left $\mathscr{A}$-module $\mathscr{A}(\,\cdot\,, v)$ is dualisable.
\end{proposition}
\noindent Proof. The following is the proof that the author learned from Geoffrey Powell. By the $\mathbb{K}$-linear version of Yoneda's lemma, we have an isomorphism of functors $\mathscr{A}\textbf{-Mod} \to \mathbf{Vect}_\mathbb{K}\colon$
	\[
		\Hom_\mathscr{A}(\mathscr{A}(\,\cdot\,, v),-) = \mathrm{Nat}(\mathscr{A}(\,\cdot\,, v),-)\cong \mathrm{ev}_v.
	\]
	Since the evaluation at $v$ is exact, so is $\Hom_\mathscr{A}(\mathscr{A}(\,\cdot\,, v),-)$; this shows that $\mathscr{A}(\,\cdot\,, v)$ is projective. Next, since any morphism $x$ in $\mathscr{A}(\,\cdot\,, v)$ is written as $x1_v$, the module $\mathscr{A}(\,\cdot\,, v)$ is generated by the sole element $1_v$. This completes the proof. \qed\\

Now everything boils down to the usual modules over an algebra by the following:
\begin{definition}
Let $\mathscr{A}$ be a $\mathbb{K}$-linear category.
\begin{itemize}
	\item The \textit{category algebra} $A$ of $\mathscr{A}$ is the (possibly non-unital) associative $\mathbb{K}$-algebra specified by the following:
	\begin{enumerate}[(i)]
		\item $A = \!\!\!\bigoplus_{v,w\in\Ob(\mathscr{A})}\!\!\! \mathscr{A}(v,w)$ as a $\mathbb{K}$-vector space with the element corresponding to $a\in\mathscr{A}$ denoted by $[a]$, and
		\item the multiplication is given, for $x\in\mathscr{A}(v_1,v_2)$ and $y\in\mathscr{A}(v_3,v_4)$, by
		\begin{align*}
			[x][y] = \left\{\begin{aligned}&[xy] &\textrm{if } v_2=v_3,\\&0& \textrm{otherwise.} \end{aligned}\right.
		\end{align*}
	\end{enumerate}
	Each $[1_v]$ is an idempotent and they are mutually orthogonal: $[1_v][1_w] = \delta_{v,w}[1_w]$. In addition, $A$ is unital if and only if $\Ob(\mathscr{A})$ is a finite set; in that case, $1_A = \sum_{v\in\Ob(\mathscr{A})}[1_v]$ gives the multiplicative unit of $A$.
	\item We set $S = S_\mathscr{A}$ to be the $\mathbb{K}$-subalgebra of $A$ generated by $[1_v]$ for all $v\in\Ob(\mathscr{A})$. If $\Ob(\mathscr{A})$ is finite, $S$ is isomorphic to the product of algebras $\prod_{\mathrm{Ob}(\mathscr{A})}\mathbb{K}$.
\end{itemize}
\end{definition}

Until the end of the section, let $\mathscr{A}$ be a $\mathbb{K}$-linear category \textit{with a finite set of objects}, and $A$ the category algebra of $\mathscr{A}$, which is unital by the assumption.

\begin{remark}\label{rem:groupoidmodule}
A left $\mathscr{A}$-module $\mathscr{N}$ is equivalent to a left $A$-module
	\[
		N =\!\!\! \bigoplus_{v\in\Ob(\mathscr{A})}\!\!\mathscr{N}(v)
	\]
	with the action given, for $x\in\mathscr{A}(u,v)$ and $n\in\mathscr{N}(w)$, by
	\begin{align*}
		[x][n] = \left\{\begin{aligned}&[xn] &\textrm{if } v=w,\\&0& \textrm{otherwise.} \end{aligned}\right.
	\end{align*}
	In the opposite direction, we can recover $\mathscr{N}$ from $N$ by the formula $\mathscr{N}(v) = [1_v]N$.
By the correspondence above, the category $\mathscr{A}\textrm{-}\mathbf{Mod}$ of left $\mathscr{A}$-modules is equivalent to the category $A\textrm{-}\mathbf{Mod}$ of left $A$-modules; in particular, the respective notion of projectivity correspond.

Similarly, an $\mathscr{A}$-bimodule is equivalent to an $A$-bimodule by replacing the role of $\mathscr{A}$ and $A$ by $\mathscr{A}^\mathrm{e}$ and $A^\mathrm{e}$, respectively.
\end{remark}

\begin{remark}\label{rem:traceproj}
Via the equivalence above, the left $\mathscr{A}$-module $\mathscr{A}(\,\cdot\,, v)$ in Proposition \ref{prop:principaldualisable} corresponds to the left $A$-module $A[1_v]$. The projectiveness of $A[1_v]$ is also directly checked as we have two left $A$-module maps
\begin{align*}
	q&\colon A\to A[1_v] \colon [a]\mapsto [a][1_v] \textrm{ and}\\
	j&\colon A[1_v] \to A \colon [a][1_v]\mapsto [a][1_v]
\end{align*}
with $q\circ j=\id$. Using these maps, the trace of $\psi\in\End_A(A[1_v])$ is given by $|(j\circ\psi\circ q)(1_A)| = |\psi([1_v])|$.\\
\end{remark}

\begin{proposition}\label{prop:grpdder}
We have an isomorphism of $\mathbb{K}$-vector spaces $\Der(\mathscr{A},\mathscr{M}) \cong \Der_S(A,M)$.
\end{proposition}
\noindent Proof. Since we set
\begin{align*}
	A = \!\!\!\bigoplus_{v,w\in\Ob(\mathscr{A})}\!\!\! \mathscr{A}(v,w)\quad \mbox{and} \quad M = \!\!\!\bigoplus_{v,w\in\Ob(\mathscr{A})}\!\!\! \mathscr{M}(v,w),
\end{align*}
a derivation $f\colon\mathscr{A}\to\mathscr{M}$ induces the $\mathbb{K}$-linear map
\[
	[f] = \!\!\!\bigoplus_{v,w\in\Ob(\mathscr{A})}\!\!\! f(v,w)\colon A\to M\colon [a]\mapsto [f(a)].
\]
We can see that $[f]$ is a derivation, and this gives the bijection we want.\qed\vspace{3pt}

\begin{corollary}
A double bracket on $\mathscr{A}$ corresponds bijectively to an $S$-linear double bracket on $A$.\qed
\end{corollary}

\begin{definition}\ 
\begin{itemize}
	\item The $\mathscr{A}$-bimodule $\Omega^1\mathscr{A}$ is defined by the universal property: $\Der(\mathscr{A},\mathscr{M})\cong \Hom_{\mathscr{A}^\mathrm{e}}(\Omega^1\mathscr{A},\mathscr{M})$ holds for all $\mathscr{M}$. By Proposition \ref{prop:grpdder}, this corresponds to $\Omega^1_SA$ through the category equivalence in Remark \ref{rem:groupoidmodule}.
	\item For a $\mathbb{K}$-linear category $\mathscr{B}$ and a left $\mathscr{B}$-module $\mathscr{M}$, a connection $\nabla\colon \mathscr{M}\to \Omega^1\mathscr{B}\otimes_\mathscr{B}\mathscr{M}$ is a collection of $\mathbb{K}$-linear maps
	\[
		\nabla(v)\colon \mathscr{M}(v)\to (\Omega^1\mathscr{B}\otimes_\mathscr{B}\mathscr{M})(v)
	\]
	satisfying the Leibniz rule. Here, the tensor product $\Omega^1\mathscr{B}\otimes_\mathscr{B}\mathscr{M}$ is defined by
	\begin{align*}
		 (\Omega^1\mathscr{B}\otimes_\mathscr{B}\mathscr{M})(v) = \bigg(\bigoplus_{w\in\mathrm{Ob}(\mathscr{B})} \Omega^1\mathscr{B}(v,w)\otimes \mathscr{M}(w)\bigg) / I
	\end{align*}
	where $I$ is the $\mathbb{K}$-vector subspace generated by
	\begin{align*}
	 \omega b\otimes m - \omega \otimes bm \;\mbox{ for }\; u,w\in\mathrm{Ob}(\mathscr{B}), \omega\in \Omega^1\mathscr{B}(v,u), b\in \mathscr{B}(u,w), m\in\mathscr{M}(w).
	\end{align*}
	Again, this is equivalent to an $S_\mathscr{B}$-linear connection
	\begin{align*}
		\nabla\colon M\to \Omega^1_{S_\mathscr{B}}B\otimes_BM.\\
	\end{align*}
\end{itemize}
\end{definition}

\section{The Loop Operation $\mu^\mathsf{fr}$}\label{sec:grpdmu}
In this section, by utilising the tools developed above, we give an algebraic description of the operation $\mu^\mathsf{fr}$ for a suitable framing $\mathsf{fr}$ on a surface (see Definition \ref{def:framing}), which will be needed to control the rotation number of a path. We first recall two maps $\kappa$ and $\mu^\mathsf{fr}$ on the groupoid algebra of the fundamental groupoid of a surface. 

Let $g,n\geq 0$ and $\Sigma = \Sigma_{g,n+1}$ be a surface of genus $g$ with $(n+1)$ boundary components $\partial_0\Sigma,\dotsc,\partial_n\Sigma$. Take a very short positive arc $\nu_j\colon [0,1] \to \partial_j\Sigma$ for each $j$ and put $\nu_j(0) = \bullet_j$, $\nu_j(1) = *_j$, and $V = \{\nu_j\}_{0\leq j\leq n}$. Now set $\mathscr{G} = \pi_1(\Sigma,V)$, the fundamental groupoid of $\Sigma$ with basepoints in $V$.

\begin{definition}
The $\mathbb{K}$-linear map $\kappa\colon \mathbb{K}\mathscr{G}\otimes \mathbb{K}\mathscr{G}\to  \mathbb{K}\mathscr{G}\otimes \mathbb{K}\mathscr{G}$ is defined as follows: for paths $\alpha\colon \bullet\to\bullet$ and $\beta\colon * \to *$ in general position,
\begin{align*}
	\kappa(\alpha,\beta) = \sum_{p\in\alpha\cap\beta} \mathrm{sign}(p;\alpha,\beta) \beta_{*p}\alpha_{p\bullet}\otimes \alpha_{\bullet p}\beta_{p*}.
\end{align*}
Here $\mathrm{sign}(p;\alpha,\beta)$ is the local intersection number ($+1$ or $-1$) with respect to the orientation of the surface.
\end{definition}

\begin{proposition}
The map $\kappa$ is a well-defined double bracket on $\mathbb{K}\mathscr{G}$.
\end{proposition}
\noindent Proof. See Section 3.2 of \cite{kappa}.\qed

\begin{figure}[tb]
\centerline{\includegraphics[scale=0.85]{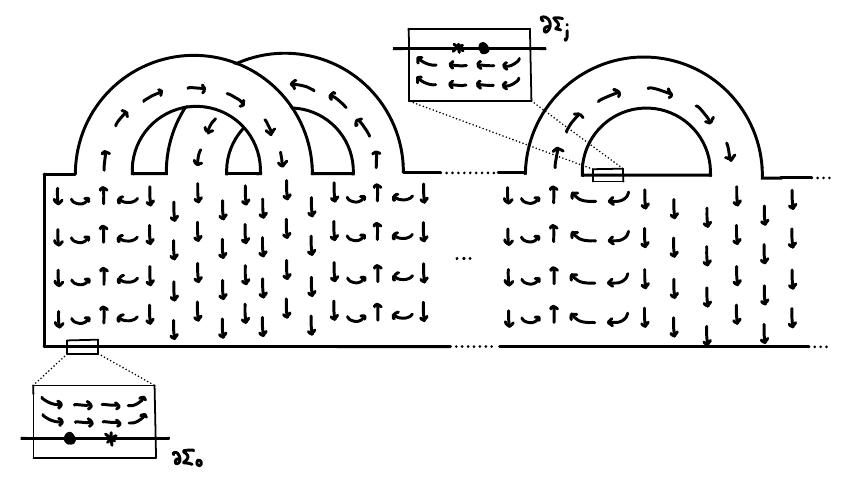}}
\centerline{\includegraphics[scale=0.85]{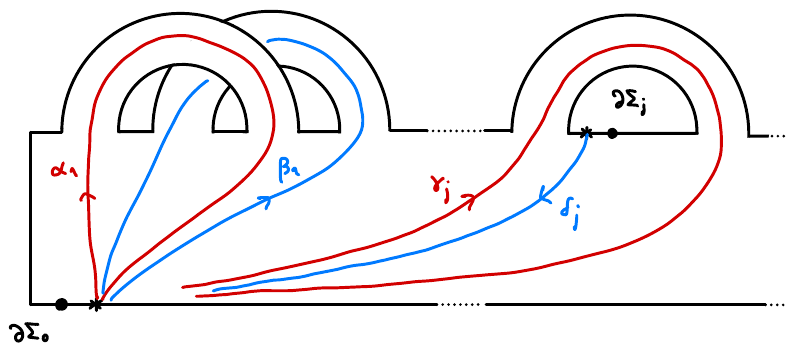}}
\caption{An example of a framing relative to $V$ (top) and the free generating system $\mathcal{C}$ (bottom).}
\label{fig:gensys}
\end{figure}

\begin{definition}\label{def:framing}
A \textit{framing} on $\Sigma$ is a nowhere-vanishing vector field. A framing is \textit{relative} to $V$ if it is positively parallel to the velocity vector of the arc $\nu_j$ in a neighbourhood of $\nu_j$ for each $j$. For a regularly immersed path $\alpha$ in $\Sigma$, we denote the rotation number with respect to $\mathsf{fr}$ by $\mathrm{rot}^\mathsf{fr}(\alpha)$.
\end{definition}

We will always take a representative of an element of $\mathscr{G}$ as a generically immersed path such that the velocity vector at each endpoint is normal to the boundary. In this case, the rotation number takes its value in $\frac12 + \mathbb{Z}$.

\begin{example}
Figure \ref{fig:gensys} (top) is an example of a framing relative to $V$. This framing is due to Alekseev--Kawazumi--Kuno--Naef; see \cite{gateder} (p.15). For the free generating system $\mathcal{C} = (\alpha_i,\beta_i,\gamma_j,\delta_j)_{1\leq i\leq g,1\leq j\leq n}$ of $\mathscr{G}$ in Figure \ref{fig:gensys} (bottom), the rotation number of a \textit{simple} path representing each generator is calculated as
\begin{align*}
	\mathrm{rot}^\mathsf{fr}(\alpha_i) = \frac12,\quad \mathrm{rot}^\mathsf{fr}(\beta_i) = -\frac12, \quad \mathrm{rot}^\mathsf{fr}(\gamma_j) = \frac12,\;\mbox{ and }\; \mathrm{rot}^\mathsf{fr}(\delta_j) = 2g+j-\frac12.
\end{align*}
\end{example}

\begin{figure}[tb]
\centerline{$+1$\hspace{70pt}$-1$}
\centerline{\includegraphics[scale=0.9]{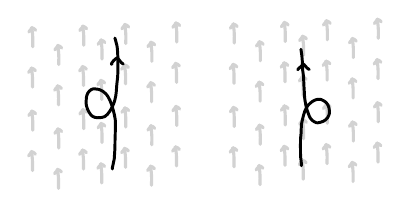}}
\caption{The left one is a positive monogon, while the right one is negative. The positive one contributes positively to the rotation number.}
\label{fig:monogon}
\end{figure}

\begin{definition}
Let $\mathsf{fr}$ be a framing on $\Sigma$ relative to $V$. The $\mathbb{K}$-linear map $\mu^\mathsf{fr}_r\colon \mathbb{K}\mathscr{G} \to |\mathbb{K}\mathscr{G}|\otimes \mathbb{K}\mathscr{G}$ is defined as follows: for a path $\alpha\in\mathscr{G}$ in general position, first deform $\alpha$ into a path from $\bullet$ to $*$ by sliding the endpoint along $\nu$, and insert positive or negative monogons so that $\mathrm{rot}^\mathsf{fr}(\alpha) = -1/2$ (see Figure \ref{fig:monogon}). Then,
\begin{align*}
	\mu^\mathsf{fr}_r(\alpha) = \sum_{p\in\mathrm{Self}(\alpha)} \mathrm{sign}(p;\alpha_\mathrm{first},\alpha_\mathrm{second}) |\alpha_{pp}|\otimes \alpha_{\bullet p*}.
\end{align*}
Here $\mathrm{Self}(\alpha)$ is the set of self-intersections of $\alpha$, and $\alpha_\mathrm{first}$ is the velocity vector of $\alpha$ passing $p$ for the first time. $\alpha_\mathrm{second}$ is analogously defined.

Similarly, the $\mathbb{K}$-linear map $\mu^\mathsf{fr}_l\colon \mathbb{K}\mathscr{G} \to \mathbb{K}\mathscr{G}\otimes |\mathbb{K}\mathscr{G}|$ is defined as follows: for a path $\alpha\in\mathscr{G}$, first deform $\alpha$ into a path from $*$ to $\bullet$ by sliding the endpoint along $\nu$, and insert positive or negative monogons so that $\mathrm{rot}^\mathsf{fr}(\alpha) = 1/2$. Then,
\begin{align*}
	\mu^\mathsf{fr}_l(\alpha) = -\sum_{p\in\mathrm{Self}(\alpha)} \mathrm{sign}(p;\alpha_\mathrm{first},\alpha_\mathrm{second}) \alpha_{*p\bullet}\otimes |\alpha_{pp}|.
\end{align*}
Set $\mu^\mathsf{fr}= \mu^\mathsf{fr}_r + \mu^\mathsf{fr}_l \colon \mathbb{K}\mathscr{G} \to |\mathbb{K}\mathscr{G}|\otimes \mathbb{K}\mathscr{G}\oplus \mathbb{K}\mathscr{G}\otimes |\mathbb{K}\mathscr{G}|$.
\end{definition}

This is essentially the same as the map $\mu$ introduced by Turaev in \cite{turaev_mu}.

\begin{proposition}\label{prop:multpropmu}
$\mu^\mathsf{fr}$ is well-defined and satisfies the following:
\[
	\mu^\mathsf{fr}(ab) = \mu^\mathsf{fr}(a)b + a\mu^\mathsf{fr}(b) + (|\cdot|\otimes \id)\kappa(a,b) + (\id\otimes|\cdot|)\kappa(b,a)\,.
\]
In particular, $\mu^\mathsf{fr}$ and $-\phi_{\kappa,\nabla}$ satisfy the same multiplicative property for any connection $\nabla$ on $\Omega^1\mathbb{K}\mathscr{G}$ by Lemma \ref{lem:multprop}.
\end{proposition}
\noindent Proof. See Section 2.3 of \cite{akkn}.\qed\\

\begin{figure}[p]
\centerline{\includegraphics[scale=0.85]{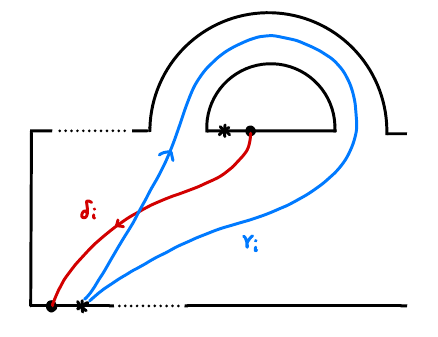} \raisebox{45pt}{$\rightsquigarrow$} \includegraphics[scale=0.85]{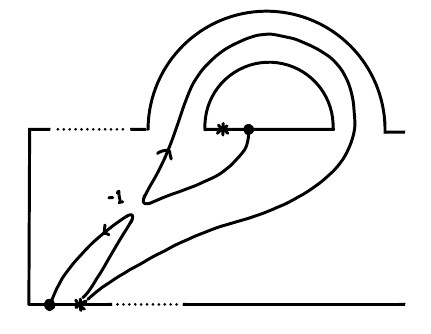}}
\caption{The computation of $\kappa(\delta_i,\gamma_i)$.}
\label{fig:i_eq_j}
\end{figure}

\begin{figure}[p]
\centerline{\phantom{$\;\rightsquigarrow$}\includegraphics[scale=0.85]{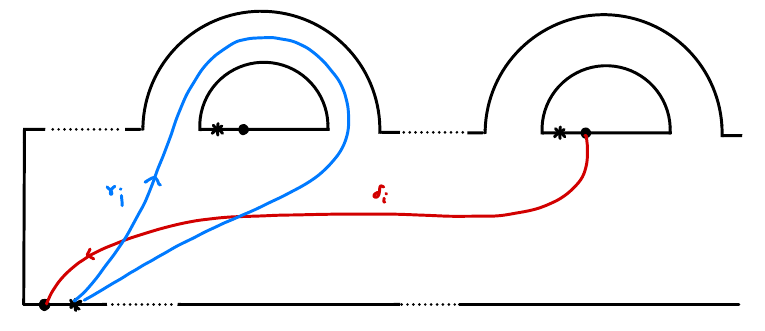}\phantom{$\rightsquigarrow$}}
\centerline{\raisebox{45pt}{$\rightsquigarrow$} \includegraphics[scale=0.85]{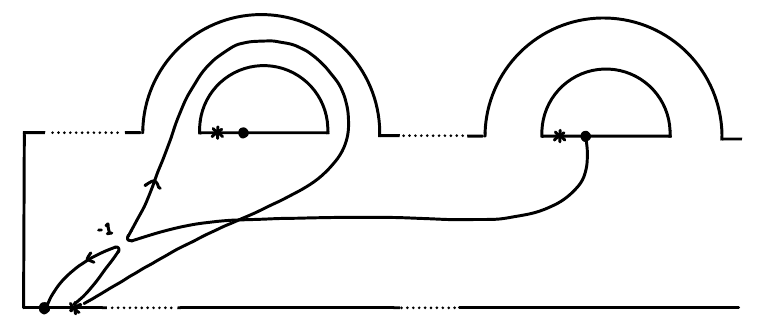}\phantom{$\;\rightsquigarrow$}}
\centerline{\raisebox{45pt}{$\hspace{40pt}+$}  \includegraphics[scale=0.85]{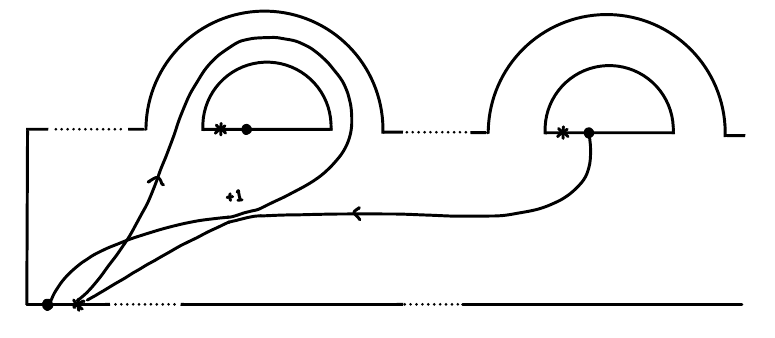}\phantom{$\;\rightsquigarrow$}}
\caption{The computation of $\kappa(\delta_i,\gamma_j)$ for $i>j$.}
\label{fig:i_geq_j}
\end{figure}

Some of the following calculations are contained in Proposition 2.16 of \cite{akkn}.
\begin{lemma}\label{lem:kappamu}
For the free generating system $\mathcal{C}$ of $\mathscr{G}$ and the framing $\mathsf{fr}$ in Figure \ref{fig:gensys}, the values of $\kappa$ and $\mu^\mathsf{fr}$ on the generators above are computed as follows:
\begin{align*}
	\kappa(\alpha_i,\alpha_j) &= \left\{\begin{aligned} &0& (i<j)\\ & \alpha_i\otimes\alpha_i - 1_0\otimes \alpha_i^2 & (i=j) \\ & \alpha_i\otimes \alpha_j + \alpha_j\otimes \alpha_i -\alpha_j\alpha_i\otimes 1_0 - 1_0\otimes \alpha_i\alpha_j& (i>j)\end{aligned}\right.,\\
	\kappa(\alpha_i,\beta_j) &= \left\{\begin{aligned} &0& (i<j)\\ & \beta_i\otimes \alpha_i& (i=j) \\ & \alpha_i\otimes \beta_j + \beta_j\otimes \alpha_i -\beta_j\alpha_i\otimes 1_0 - 1_0\otimes \alpha_i\beta_j& (i>j)\end{aligned}\right.,\\
	\kappa(\alpha_i,\gamma_j) &= 0,\\	
	\kappa(\alpha_i,\delta_j) &= 0,\\
	\kappa(\beta_i,\alpha_j) &= \left\{\begin{aligned} &0& (i<j)\\ & \beta_i\otimes\alpha_i - \alpha_i\beta_i\otimes 1_0 - 1_0\otimes\beta_i\alpha_i & (i=j) \\ & \beta_i\otimes \alpha_j + \alpha_j\otimes \beta_i -\alpha_j\beta_i\otimes 1_0 - 1_0\otimes \beta_i\alpha_j& (i>j)\end{aligned}\right.,\\
	\kappa(\beta_i,\beta_j) &= \left\{\begin{aligned} &0& (i<j)\\ & \beta_i\otimes\beta_i - \beta_i^2\otimes 1_0& (i=j) \\ & \beta_i\otimes \beta_j + \beta_j\otimes \beta_i -\beta_j\beta_i\otimes 1_0 - 1_0\otimes \beta_i\beta_j& (i>j)\end{aligned}\right.,\\
	\kappa(\beta_i,\gamma_j) &= 0,\\
	\kappa(\beta_i,\delta_j) &= 0,\\
	\kappa(\gamma_i,\alpha_j) &= \gamma_i\otimes \alpha_j + \alpha_j\otimes \gamma_i -\alpha_j\gamma_i\otimes 1_0 - 1_0\otimes \gamma_i\alpha_j,\\
	\kappa(\gamma_i,\beta_j) &= \gamma_i\otimes \beta_j + \beta_j\otimes \gamma_i -\beta_j\gamma_i\otimes 1_0 - 1_0\otimes \gamma_i\beta_j,\\
	\kappa(\gamma_i,\gamma_j) &= \left\{\begin{aligned} &0& (i<j)\\ & \gamma_i\otimes\gamma_i -1_0\otimes \gamma_i^2& (i=j) \\ & \gamma_i\otimes \gamma_j + \gamma_j\otimes \gamma_i -\gamma_j\gamma_i\otimes 1_0 - 1_0\otimes \gamma_i\gamma_j& (i>j)\end{aligned}\right.,\\
	\kappa(\gamma_i,\delta_j) &= \left\{\begin{aligned} &0& (i<j)\\ & \delta_i\otimes \gamma_i & (i=j) \\ & -\delta_j\gamma_i\otimes 1_0 + \delta_j\otimes \gamma_i & (i>j)\end{aligned}\right.,\\
	\kappa(\delta_i,\alpha_j) &= \alpha_j\otimes \delta_i - 1_0\otimes \delta_i\alpha_j,\\
	\kappa(\delta_i,\beta_j) &= \beta_j\otimes \delta_i - 1_0\otimes \delta_i\beta_j,\\
	\kappa(\delta_i,\gamma_j) &= \left\{\begin{aligned} &0& (i<j)\\ &-1_0\otimes \delta_i\gamma_i & (i=j) \\ &\gamma_j\otimes \delta_i - 1_0\otimes \delta_i\gamma_j& (i>j)\end{aligned}\right.,\\
	\kappa(\delta_i,\delta_j) &= \left\{\begin{aligned} &0& (i<j)\\ & \delta_i\otimes \delta_i & (i=j) \\ & \delta_j\otimes \delta_i & (i>j)\end{aligned}\right.,\\
	\mu^\mathsf{fr}(\alpha_i) &= |1_0|\otimes \alpha_i - 1_0\otimes |\alpha_i|,\\
	\mu^\mathsf{fr}(\beta_i) &= - |\beta_i| \otimes 1_0 + \beta_i\otimes |1_0|,\\
	\mu^\mathsf{fr}(\gamma_i) &= |1_0|\otimes \gamma_i - 1_0\otimes |\gamma_i|,\\
	\mu^\mathsf{fr}(\delta_i) &= (2g+i) |1_0|\otimes \delta_i - (2g+i-1)\delta_i \otimes |1_0|.
\end{align*} 
\end{lemma}
\noindent Proof. We only compute $\kappa(\delta_i,\gamma_j)$ and $\mu^\mathsf{fr}(\delta_i)$; the other cases can be similarly done. 

For $\kappa(\delta_i,\gamma_j)$, take simple representatives of $\delta_i$ and $\gamma_j$ as in Figure \ref{fig:gensys}. If $i<j$, $\delta_i$ and $\gamma_j$ are disjoint, therefore $\kappa(\delta_i,\gamma_j) = 0$. If $i=j$, we compute as in Figure \ref{fig:i_eq_j}; the local intersection number at the only intersection is $-1$ with respect to the counterclockwise orientation, and each loop reads $1$ and $\delta_i\gamma_i$, respectively. If $i>j$, we compute as in Figure \ref{fig:i_geq_j}; the rest is similar.

For $\mu^\mathsf{fr}(\delta_i)$, we first compute $\mu_r^\mathsf{fr}(\delta_i)$. Since $\mathrm{rot}^\mathsf{fr}(\delta_i) = 2g+i-\frac12$, we insert $(2g+i)$ negative monogons to make the rotation number $-\frac12$. Each negative monogon contributes $|1_0|\otimes \delta_i$ with the coefficient $+1$; note that we have $|1_0| = |1_i|$ for any $i$. Next, $\mu_l^\mathsf{fr}(\delta_i)$ is computed by inserting $(2g+i-1)$ negative monogons so that $\mathrm{rot}^\mathsf{fr} = \frac12$. Since $\mu_l^\mathsf{fr}$ has got a $(-1)$ in its definition, each negative monogon contributes $\delta_i\otimes |1_0|$ with the coefficient $(-1)$. This completes the proof. \qed\\

\begin{lemma}
$\Omega^1\mathbb{K}\mathscr{G}$ is dualisable over $\mathbb{K}\mathscr{G}^\mathrm{e}$.
\end{lemma}
\noindent Proof. For $c\in \mathcal{C}$, we define a double derivation $\partial_c\colon \mathbb{K}\mathscr{G} \to \mathbb{K}\mathscr{G}\otimes \mathbb{K}\mathscr{G}$ by
\[
	\partial_c(c') = \delta_{cc'}\cdot 1_{s(c')}\otimes 1_{t(c')}\;\mbox{ for }\; c'\in\mathcal{C},
\]
where $\delta_{cc'}$ is Kronecker's delta, $s(c')$ is the source of $c'$ and $t(c')$ is the target of $c'$. This is well-defined since $\mathscr{G}$ is the free groupoid generated by $\mathcal{C}$. For any $\omega\in\Omega^1\mathbb{K}\mathscr{G}$, we have
\[
	\omega = \sum_{c\in\mathcal{C}} i_{\partial_c} (\omega)' \, dc \, i_{\partial_c} (\omega)''
\]
Therefore, we have the following isomorphism of $\mathbb{K}\mathscr{G}$-bimodules:
\begin{align*}
	\Omega^1\mathbb{K}\mathscr{G} \cong \bigoplus_{c = \alpha_i,\beta_i,\gamma_j} \mathbb{K}\mathscr{G}(\argdot,\bullet_0)\,(dc)c^{-1}\, \mathbb{K}\mathscr{G}(\bullet_0,\argdot) \oplus \bigoplus_{1\leq j\leq n} \mathbb{K}\mathscr{G}(\argdot,\bullet_j)\,(d\delta_j)\delta_j^{-1} \,\mathbb{K}\mathscr{G}(\bullet_j,\argdot)\,.
\end{align*}
Each summand is isomorphic to the left $\mathbb{K}\mathscr{G}^\mathrm{e}$-module of the form $\mathbb{K}\mathscr{G}^\mathrm{e}(\,\cdot\,, (\bullet_i,\bullet_i))$, which is dualisable by Lemma \ref{prop:principaldualisable}.\qed\\

\begin{remark}\label{rem:omegaKGtrace}
The trace of $\psi\in\End_{\mathbb{K}\mathscr{G}^\mathrm{e}}(\Omega^1\mathbb{K}\mathscr{G})$ is given by the sum of the traces on each summand. Therefore, by Remark \ref{rem:traceproj}, we obtain
\[
	\Tr(\psi) = \sum_{c \in\mathcal{C}} \left| (\textrm{the coefficient of } (dc)c^{-1}\textrm{ in }\psi((dc)c^{-1}) ) \right|\,.\\
\]
\end{remark}

The following is the main result of this paper.

\begin{theorem}\label{thm:groupoidmu}
Let the free-generating system $\mathcal{C}$ of $\mathscr{G}$ and the framing $\mathsf{fr}$ normalised near the base points $V$ be the ones given in Figure \ref{fig:gensys}, and define the connection $\nabla\!_\mathcal{C}$ on $\Omega^1\mathbb{K}\mathscr{G}$ by $\nabla\!_\mathcal{C}((dc)c^{-1}) = 0$ for all $c\in\mathcal{C}$. Then the map $-\phi_{\kappa,\nabla\!_\mathcal{C}}$ defined in \ref{def:phipinabla} coincides with the framed version $\mu^\mathsf{fr}$ of Turaev's loop operation.
\end{theorem}
\noindent Proof. Recall that $-\phi_{\kappa,\nabla\!_\mathcal{C}}$ is the trace of $L_{\kappa(c)} - i_{\kappa(c)^\mathrm{e}}\nabla\!_\mathcal{C}$. We compute the values of $L_{\kappa(c)}((dc')c'^{-1})$ in $\Omega^1\mathbb{K}\mathscr{G}$ for $c,c'\in\mathcal{C}$. First of all, we have
\[
	L_{\kappa(c)}((dc')c'^{-1}) = d\kappa(c,c')c'^{-1} - (dc')c'^{-1}\kappa(c,c')c'^{-1}\,.
\]
By Lemma \ref{lem:kappamu}, we have,
\begin{align*}
	L_{\kappa(\alpha_i)}((d\alpha_j)\alpha_j^{-1}) &= \left\{\begin{aligned}
		&0&(i<j)\phantom{,} \\
		&(\alpha_i\otimes 1_0\otimes \bar 1_0 - 1_0\otimes 1_0\otimes \bar \alpha_i -1_0\otimes \alpha_i\otimes \bar 1_0 + 1_0\otimes \bar \alpha_i\otimes \bar 1_0)\cdot (d\alpha_i)\alpha_i^{-1} & (i=j)\phantom{,}\\
		&\begin{aligned}
			&( \alpha_i\otimes 1_0\otimes \bar 1_0 - 1_0\otimes \alpha_i\otimes \bar 1_0 - 1_0\otimes \bar 1_0\otimes \bar \alpha_i + 1_0\otimes \bar \alpha_i\otimes \bar 1_0 )\cdot (d\alpha_j)\alpha_j^{-1}\\
			&\qquad + (\cdots)\cdot (d\alpha_i)\alpha_i^{-1} \end{aligned}& (i>j),
	 \end{aligned}\right.\\
	 L_{\kappa(\alpha_i)}((d\beta_j)\beta_j^{-1}) &= \left\{\begin{aligned}
		&0&(i<j)\phantom{,} \\
		& (\cdots)\cdot (d\alpha_i)\alpha_i^{-1} & (i=j)\phantom{,}\\
		&\begin{aligned}
			&( \alpha_i\otimes 1_0\otimes \bar 1_0 - 1_0\otimes \alpha_i\otimes \bar 1_0 - 1_0\otimes \bar 1_0\otimes \bar \alpha_i + 1_0\otimes \bar \alpha_i\otimes \bar 1_0 )\cdot (d\beta_j)\beta_j^{-1}\\
			&\qquad + (\cdots)\cdot (d\alpha_i)\alpha_i^{-1}\end{aligned} & (i>j),
	 \end{aligned}\right.\\
	 L_{\kappa(\alpha_i)}((d\gamma_j)\gamma_j^{-1}) &= 0,\\
	 L_{\kappa(\alpha_i)}((d\delta_j)\delta_j^{-1}) &= 0.
\end{align*}
For example, $L_{\kappa(\alpha_i)}((d\alpha_i)\alpha_i^{-1})$ is computed as follows:
\begin{align*}
	L_{\kappa(\alpha_i)}((d\alpha_i)\alpha_i^{-1}) &= d\kappa(\alpha_i,\alpha_i)\alpha_i^{-1} - (d\alpha_i)\alpha_i^{-1}\kappa(\alpha_i,\alpha_i)\alpha_i^{-1}\\
	&= d(\alpha_i\otimes\alpha_i - 1_0\otimes \alpha_i^2)\alpha_i^{-1} - (d\alpha_i)\alpha_i^{-1}(\alpha_i\otimes\alpha_i - 1_0\otimes \alpha_i^2)\alpha_i^{-1}\\
	&= (d\alpha_i\otimes\alpha_i + \alpha_i\otimes d\alpha_i  - 1_0\otimes (d\alpha_i)\alpha_i - 1_0\otimes \alpha_i d\alpha_i)\alpha_i^{-1} - (d\alpha_i)\alpha_i^{-1}(\alpha_i\otimes\alpha_i - 1_0\otimes \alpha_i^2)\alpha_i^{-1}\\
	&= (d\alpha_i\otimes 1_0+ \alpha_i\otimes (d\alpha_i)\alpha_i^{-1} - 1_0\otimes d\alpha_i - 1_0\otimes \alpha_i (d\alpha_i)\alpha_i^{-1}) -( d\alpha_i \otimes 1_0 - (d\alpha_i)\alpha_i^{-1}\otimes \alpha_i)\\
	&= \alpha_i\otimes (d\alpha_i)\alpha_i^{-1} - 1_0\otimes d\alpha_i - 1_0\otimes \alpha_i (d\alpha_i)\alpha_i^{-1} +  (d\alpha_i)\alpha_i^{-1}\otimes \alpha_i\\
	&= (\alpha_i\otimes 1_0\otimes \bar 1_0 - 1_0\otimes 1_0\otimes \bar \alpha_i - 1_0\otimes \alpha_i\otimes \bar 1_0 + 1_0\otimes \bar \alpha_i \otimes \bar 1_0)\cdot (d\alpha_i)\alpha_i^{-1},
\end{align*}
where the last equality comes from the identification in Lemma \ref{lem:tripM}. Since $-\phi_{\kappa,\nabla\!_\mathcal{C}}(\alpha_i)$ is the trace of $(L_{\kappa(\alpha_i)} - i_{\kappa(\alpha_i)^\mathrm{e}}\nabla\!_\mathcal{C})$, it is given by the formula in Remark \ref{rem:omegaKGtrace}. For the term involving $L_{\kappa(\alpha_i)}((d\alpha_i)\alpha_i^{-1})$, the coefficient is computed as
\begin{align*}
	|\alpha_i\otimes 1_0\otimes \bar 1_0 - 1_0\otimes 1_0\otimes \bar \alpha_i -1_0\otimes \alpha_i\otimes \bar 1_0 + 1_0\otimes \bar \alpha_i\otimes \bar 1_0| &= \alpha_i \otimes |1_0| - 1_0\otimes |\alpha_i| - \alpha_i \otimes |1_0| + |1_0|\otimes \alpha_i\\
	&= - 1_0\otimes |\alpha_i|  + |1_0|\otimes \alpha_i
\end{align*}
via the isomorphism (\ref{eq:tripA}). The other coefficients vanish upon taking the cyclic quotient. For example, for $L_{\kappa(\alpha_i)}((d\alpha_j)\alpha_j^{-1})$, $i\neq j$, we have
\begin{align*}
	&|\alpha_i\otimes 1_0\otimes \bar 1_0 - 1_0\otimes \alpha_i\otimes \bar 1_0 - 1_0\otimes \bar 1_0\otimes \bar \alpha_i + 1_0\otimes \bar \alpha_i\otimes \bar 1_0|\\
	&\quad = \alpha_i\otimes |1_0| - \alpha_i\otimes |1_0| - |1_0|\otimes \alpha_i + |1_0|\otimes \alpha_i = 0\,.
\end{align*}
Therefore, we obtain $-\phi_{\kappa,\nabla\!_\mathcal{C}}(\alpha_i) = - 1_0\otimes |\alpha_i| + |1_0|\otimes \alpha_i = \mu^\mathsf{fr}(\alpha_i)$ by Lemma \ref{lem:kappamu}. The rest is similar, so we just show the result of the computation. We have
\begin{align*}
	L_{\kappa(\beta_i)}((d\alpha_j)\alpha_j^{-1}) &= \left\{\begin{aligned}
		&0&(i<j)\phantom{,} \\
		&\begin{aligned}
			&( \beta_i\otimes 1_0\otimes \bar 1_0 - 1_0\otimes \beta_i\otimes \bar 1_0 - 1_0\otimes \bar 1_0\otimes \bar \beta_i + 1_0\otimes \bar \beta_i\otimes \bar 1_0 )\cdot (d\alpha_j)\alpha_j^{-1}\\
			&\qquad + (\cdots)\cdot (d\beta_i)\beta_i^{-1}\end{aligned} & (i\geq j),
	\end{aligned}\right.\\
	L_{\kappa(\beta_i)}((d\beta_j)\beta_j^{-1}) &= \left\{\begin{aligned}
		&0&(i<j)\phantom{,} \\
		& ( \beta_i\otimes 1_0\otimes \bar 1_0 - \beta_i\otimes \bar\beta_i^{-1}\otimes \bar\beta_i )\cdot (d\beta_i)\beta_i^{-1}& (i=j)\phantom{,}\\
		&\begin{aligned}
			&( \beta_i\otimes 1_0\otimes \bar 1_0 - 1_0\otimes \beta_i\otimes \bar 1_0 - 1_0\otimes \bar 1_0\otimes \bar \beta_i + 1_0\otimes \bar \beta_i\otimes \bar 1_0 )\cdot (d\beta_j)\beta_j^{-1}\\
			&\qquad + (\cdots)\cdot (d\beta_i)\beta_i^{-1} \end{aligned}& (i>j),
	\end{aligned}\right.\\
	L_{\kappa(\beta_i)}((d\gamma_j)\gamma_j^{-1}) &= 0,\\
	L_{\kappa(\beta_i)}((d\delta_j)\delta_j^{-1}) &= 0.
\end{align*}
Therefore, $-\phi_{\kappa,\nabla\!_\mathcal{C}}(\beta_i) = \beta_i\otimes |1_0| - |\beta_i|\otimes 1_0 = \mu^\mathsf{fr}(\beta_i)$. Next,
\begin{align*}
	L_{\kappa(\gamma_i)}((d\alpha_j)\alpha_j^{-1}) &= ( \gamma_i\otimes 1_0\otimes \bar 1_0 - 1_0\otimes \gamma_i\otimes \bar 1_0 - 1_0\otimes \bar 1_0\otimes \bar \gamma_i + 1_0\otimes \bar \gamma_i\otimes \bar 1_0 )\cdot (d\alpha_j)\alpha_j^{-1} + (\cdots)\cdot (d\gamma_i)\gamma_i^{-1},\\
	L_{\kappa(\gamma_i)}((d\beta_j)\beta_j^{-1}) &= ( \gamma_i\otimes 1_0\otimes \bar 1_0 - 1_0\otimes \gamma_i\otimes \bar 1_0 - 1_0\otimes \bar 1_0\otimes \bar \gamma_i + 1_0\otimes \bar \gamma_i\otimes \bar 1_0 )\cdot (d\beta_j)\beta_j^{-1} + (\cdots)\cdot (d\gamma_i)\gamma_i^{-1},\\
	L_{\kappa(\gamma_i)}((d\gamma_j)\gamma_j^{-1}) &= \left\{\begin{aligned}
		&0&(i<j)\phantom{,} \\
		& ( \gamma_i\otimes 1_0\otimes \bar 1_0 - 1_0\otimes 1_0\otimes \bar\gamma_i - 1_0\otimes \gamma_i\otimes \bar 1_0 + 1_0\otimes \bar\gamma_i\otimes\bar 1_0)\cdot (d\gamma_i)\gamma_i^{-1}& (i=j)\phantom{,}\\
		&\begin{aligned}
			&( \gamma_i\otimes 1_0\otimes \bar 1_0 - 1_0\otimes \gamma_i\otimes \bar 1_0 - 1_0\otimes \bar 1_0\otimes \bar \gamma_i + 1_0\otimes \bar \gamma_i\otimes \bar 1_0 )\cdot (d\gamma_j)\gamma_j^{-1}\\
			&\qquad + (\cdots)\cdot (d\gamma_i)\gamma_i^{-1}\end{aligned} & (i>j),
	\end{aligned}\right.\\
	L_{\kappa(\gamma_i)}((d\delta_j)\delta_j^{-1}) &= \left\{\begin{aligned}
		&0&(i<j) \\
		& (\cdots)\cdot (d\gamma_i)\gamma_i^{-1}& (i=j)\\
		& (\cdots)\cdot (d\gamma_i)\gamma_i^{-1}& (i>j)
	\end{aligned}\right..
\end{align*}
Therefore, $-\phi_{\kappa,\nabla\!_\mathcal{C}}(\gamma_i) = -1_0\otimes |\gamma_i| + |1_0|\otimes \gamma_i =  \mu^\mathsf{fr}(\gamma_i)$. Finally,
\begin{align*}
	L_{\kappa(\delta_i)}((d\alpha_j)\alpha_j^{-1}) &= (-1_0\otimes \delta_i\otimes \bar 1_0 + 1_0\otimes \bar \delta_i\otimes \bar 1_0)\cdot (d\alpha_j)\alpha_j^{-1} + (\cdots)\cdot (d\delta_i)\delta_i^{-1},\\
	L_{\kappa(\delta_i)}((d\beta_j)\beta_j^{-1}) &= (-1_0\otimes \delta_i\otimes \bar 1_0 + 1_0\otimes \bar \delta_i\otimes \bar 1_0)\cdot (d\beta_j)\beta_j^{-1} + (\cdots)\cdot (d\delta_i)\delta_i^{-1},\\
	L_{\kappa(\delta_i)}((d\gamma_j)\gamma_j^{-1}) &= \left\{\begin{aligned}
		&0&(i<j) \\
		& (-1_0\otimes \delta_i\otimes \bar 1_0 + 1_0\otimes \bar \delta_i\otimes \bar 1_0)\cdot (d\gamma_i)\gamma_i^{-1} &(i=j)\\
		&(-1_0\otimes \delta_i\otimes \bar 1_0 + 1_0\otimes \bar \delta_i\otimes \bar 1_0)\cdot (d\gamma_j)\gamma_j^{-1} + (\cdots)\cdot (d\delta_i)\delta_i^{-1} & (i>j)
	\end{aligned}\right.,\\
	L_{\kappa(\delta_i)}((d\delta_j)\delta_j^{-1}) &= \left\{\begin{aligned}
		&0&(i<j) \\
		& (\delta_i\otimes 1_i\otimes \bar 1_i)\cdot d\delta_i \delta_i^{-1} & (i=j)\\
		& 0 & (i>j)
	\end{aligned}\right..
\end{align*}
Therefore, using $|1_0| = |1_i|$, we have
\begin{align*}
	-\phi_{\kappa,\nabla\!_\mathcal{C}}(\delta_i) &= (2g+i)(-\delta_i\otimes |1_0| + |1_0|\otimes \delta_i) + \delta_i\otimes |1_i|\\
	&= (2g+i)|1_0|\otimes \delta_i - (2g + i - 1)\delta_i\otimes |1_0| = \mu^\mathsf{fr}(\delta_i).
\end{align*}
Since $\mu^\mathsf{fr}$ and $-\phi_{\kappa,\nabla\!_\mathcal{C}}$ satisfy the same multiplicative property by Lemma \ref{lem:multprop} and Proposition \ref{prop:multpropmu}, they coincide on the whole $\mathbb{K}\mathscr{G}$.\qed\\

Now we calculate the modular vector field associated with $\kappa$ and $\nabla\!_\mathcal{C}$.
\begin{lemma}
The double bracket $\kappa$ satisfies $\kappa + \kappa^\circ = \ad_{\ad_e}$ with $e=\sum_{0\leq j\leq n}[1_j]\otimes [1_j]$. 
\end{lemma}
\noindent Proof. Since both $\kappa + \kappa^\circ$ and $ \ad_{\ad_e}$ are double derivations, we only have to check that they agree on the generating system $\mathcal{C}$. First of all, we have
\[
	\ad_{\ad_e}(x,y) = -ye'x\otimes e'' + e'x\otimes e''y + ye'\otimes xe'' - e'\otimes xe''y.
\]
for $x,y\in A$. If $x,y\in\{[\alpha_i],[\beta_i],[\gamma_j]\}$, their endpoints are always $1_0$ so only the terms involving $1_0$ survive:\begin{align*}
	\mathrm{(RHS)} &= -y[1_0]x\otimes [1_0] + [1_0]x\otimes [1_0]y + y[1_0]\otimes x[1_0] - [1_0]\otimes x[1_0]y\\
	&= -yx\otimes [1_0] + x\otimes y + y\otimes x - [1_0]\otimes xy\\
	&= (\kappa + \kappa^\circ)(x,y)
\end{align*}
by Lemma \ref{lem:kappamu}. Next, if $x\in\{[\alpha_i],[\beta_i],[\gamma_j]\}$ and $y=[\delta_k]$, we have
\begin{align*}
	\mathrm{(RHS)} &= -[\delta_k][1_0]x\otimes [1_0] + 0 + [\delta_k][1_0]\otimes x[1_0] - 0\\
	&= -[\delta_k]x\otimes [1_0] + [\delta_k]\otimes x\\
	&= (\kappa + \kappa^\circ)(x,[\delta_k]).
\end{align*}
Finally, if $x = [\delta_i]$ and $y = [\delta_j]$, we have
\begin{align*}
	\mathrm{(RHS)} &= [\delta_j][1_0]\otimes [\delta_i][1_0] + \sum_k [1_k][\delta_i]\otimes [1_k][\delta_j].
\end{align*}
If $i>j$, this is equal to $[\delta_j]\otimes [\delta_i]$, while if $i=j$, this is equal to $2[\delta_i]\otimes[\delta_i]$. In either way, this is equal to $(\kappa + \kappa^\circ)([\delta_i],[\delta_j])$. Since both $\kappa + \kappa^\circ$ and $ \ad_{\ad_e}$ are symmetric, we have checked all the necessary pairs. \qed

\begin{theorem}\label{thm:mkappanabla}
We have $\mathbf{m}_{\kappa,\nabla\!_\mathcal{C},\ad_e} = 0$.
\end{theorem}
\noindent Proof. If $c\in\{\alpha_i,\beta_i,\gamma_j\}$, by Theorem \ref{thm:groupoidmu} together with Lemma \ref{lem:kappamu}, we have
\begin{align*}
	\mathrm{fd}(\phi_{\kappa,\nabla\!_\mathcal{C}}([c])) = [1_0]\otimes |[c]| - [c]\otimes |[1_0]|
\end{align*}
and therefore
\begin{align*}
	\mathbf{m}_{\kappa,\nabla\!_\mathcal{C},\ad_e}([c]) &= \mathrm{fd}(\phi_{\kappa,\nabla\!_\mathcal{C}}([c])) + (\id\otimes|\cdot|)\circ\ad_e([c])\\
	&= ([1_0]\otimes |[c]| - c\otimes |[1_0]|) + ([c]\otimes |[1_0]| - [1_0]\otimes |[c]|) = 0.
\end{align*}
If $c = \delta_j$, again  by Theorem \ref{thm:groupoidmu} together with Lemma \ref{lem:kappamu}, we have
\[
	\mathrm{fd}(\phi_{\kappa,\nabla\!_\mathcal{C}}([\delta_j])) = -[\delta_j]\otimes |[1_0]|
\]
and
\begin{align*}
	\mathbf{m}_{\kappa,\nabla\!_\mathcal{C},\ad_e}([\delta_j]) &= \mathrm{fd}(\phi_{\kappa,\nabla\!_\mathcal{C}}([\delta_j])) + (\id\otimes|\cdot|)\circ\ad_e([\delta_j])\\
	&= -[\delta_j]\otimes |[1_0]| + ([\delta_j]\otimes |[1_0]| - [1_j]\otimes |[\delta_j]|).
\end{align*}
In addition, $|[\delta_j]|$ in the last term is zero since we have $|[\delta_j]| = |[1_j][\delta_j]| = |[\delta_j][1_j]| = 0$ in $|A|$. Therefore, we have $\mathbf{m}_{\kappa,\nabla\!_\mathcal{C},\ad_e}([\delta_j]) = 0$. Finally, since $\mathbf{m}_{\kappa,\nabla\!_\mathcal{C},\ad_e}$ is a derivation by Definition-Lemma \ref{deflem:modvec}, $\mathbf{m}_{\kappa,\nabla\!_\mathcal{C},\ad_e}$ vanishes identically on the category algebra of $\mathbb{K}\mathscr{G}$. \qed

\begin{corollary}
Putting $\sigma = \mathrm{Ham}_\kappa$, the map $\delta^{\sigma,\nabla\!_\mathcal{C}}$ is skew-symmetric.
\end{corollary}
\noindent Proof. This follows from Theorem \ref{thm:modvec} together with Theorem \ref{thm:mkappanabla}.\qed

\begin{remark}
By Theorem \ref{thm:groupoidmu}, we have the following equality
\[
	\delta^{\sigma,\nabla\!_\mathcal{C}} = |\phi_{\kappa,\nabla\!_\mathcal{C}}| = -|\mu^\mathsf{fr}| = -\delta_\mathrm{T}^\mathsf{fr},
\]
where $\delta_\mathrm{T}^\mathsf{fr}$ is the framed version of the Turaev cobracket; see Equation (12) of \cite{akkn} for the definition. Then, the corollary above can also be deduced from the skew-symmetry of the Turaev cobracket. 
\end{remark}

\small
\bibliographystyle{alphaurl}
\bibliography{tm.bib}

\end{document}